\documentclass[amscd,amssymb,verbatim,12pt]{amsart}
\usepackage[all]{xy}
\usepackage{amssymb,latexsym, amsmath, amscd, array, graphicx}
\xyoption{v2}
\newtheorem{dummy}{realdumb}[section]
\newtheorem{thm}[dummy]{Theorem}
\newtheorem{lem}[dummy]{Lemma}
\newtheorem{cor}[dummy]{Corollary}
\newtheorem{prop}[dummy]{Proposition}

\theoremstyle{definition}       
\newtheorem{defn}[dummy]{Definition}

\newtheorem{ques}[dummy]{Question}

\newtheorem{rem}[dummy]{Remark}

\newtheoremstyle{break}
   {9pt}
   {9pt}
   {\rmfamily}
   {}
   {\bfseries}
   {.}
   {\newline}
{}

\theoremstyle{break}

\newtheorem{brem}[dummy]{Remark}

\DeclareMathOperator{\id}{id}
\DeclareMathOperator{\T}{T}
\DeclareMathOperator{\im}{im}
\DeclareMathOperator{\rel}{rel}
\DeclareMathOperator{\Ker}{Ker}
\DeclareMathOperator{\Emb}{Emb}
\DeclareMathOperator{\GTop}{G/TOP}
\DeclareMathOperator{\dist}{dist}
\DeclareMathOperator{\Tor}{Tor}

\newcommand{\CP}{\mathbb {CP}}

\def\ni{\noindent}

\def\bQ{\mathbb Q}
\def\bR{\mathbb R}
\def\bE{\mathbb E}
\def\bH{\mathbb H}
\def\bL{\mathbb L}
\def\bC{\mathbb C}
\def\bP{\mathbb P}
\def\ni{\noindent}
\def\x{\times}

\def\mk{\vskip .1truein}

\def\p{\partial}

\def\KO{\overline{KO}_*}

\def\SS(#1,#2){\mathcal S{\bigl(\begin{smallmatrix} #1\\ \downarrow \\ #2
\end{smallmatrix} \bigr)}}

\def\cm{\mathcal C\mathcal M}
\def\mn{{\mathcal M^{man}}(n,\rho)}
\def\bZ{\mathbb Z}
\def\circM{\overset \circ M}
\def\circD{\overset \circ D}
\def\CS{\mathcal S}
\def\CB{\mathcal B}

\begin{document}

\title{An infinite-dimensional phenomenon in finite-dimensional metric topology}

\author {Alexander N. Dranishnikov}
\address {Department of Mathematics, University of Florida,
Gainesville, FL 32611-8105}

\author{Steven C. Ferry}
\address {Department of Mathematics, Rutgers University, 
Piscataway, NJ 08854-8019}

\author{Shmuel Weinberger}
\address {Department of Mathematics, University of Chicago, Chicago, IL 60737}

\thanks{The first and third authors are partially supported by NSF grants.  The first two authors
would like to thank Max Planck Institute for Mathematics for hospitality and
excellent working conditions. The second
author would like to thank the University of
Chicago for hospitality during numerous visits.}

\date{\today}

\begin{abstract}
We show that there are homotopy equivalences $h:N\to M$ between closed
manifolds which are induced by cell-like maps $p:N\to X$ and $q:M\to
X$ but which are not homotopic to homeomorphisms.  The phenomenon is
based on construction of cell-like maps that kill certain
$\bL$-classes.  The image space in these constructions is necessarily infinite-dimensional. In dimension $>6$ we classify all such homotopy equivalences. 
As an application, we show that such homotopy equivalences are realized by 
deformations of Riemannian manifolds in Gromov-Hausdorff space preserving 
a contractibility function.
\end{abstract}

\maketitle

\section{Introduction}

\ni The authors were led to the questions studied in this paper by two different routes.  The first route was via the theory of CE-maps, homology manifolds, and resolutions and the second route was via a quest to understand precompact subsets of Riemannian manifolds in Gromov-Hausdorff space.  Connecting these problems to each other led to a new functorial subgroup of the structure group of topological surgery theory and to examples casting light on  both of these problems.
\mk
\ni Beginning with the second question, recall that the Gromov-Hausdorff metric is a metric on the isomorphism classes of compact metric spaces.  The Gromov-Hausdorff distance from a metric space X to the one-point metric space p is diam(X)/2, so being Gromov-Hausdorff close imposes little connection between the topologies of compact metric spaces.
\mk
\ni
However, if one assumes a uniform local contractibility condition, then much more structure is preserved.  Let $\rho:[0,R) \to [0, \infty)$ be a function with $\rho(0)=0$ and $\rho(t)\ge t$, such that $\rho$ is continuous at 0.  Following Borsuk \cite{B} and Gromov \cite{Gromov}, we say that $X$ is $LGC(\rho)$ if every ball of radius $r < R$ in $X$ is nullhomotopic in the concentric ball of radius $\rho(r)$.  This is a generalization of the idea of injectivity radius for Riemannian manifolds.  Sufficiently Gromov-Hausdorff close n-dimensional $LGC(\rho)$ spaces are homotopy equivalent -- and there are explicit estimates on the required degree of closeness in terms of $n$ and $\rho$. See \cite{B}, \cite{Pet}.
\mk
\ni 	A theorem of Chapman and Ferry \cite{ChF} implies that if $M$ is a closed metric $n$-manifold, $n\ge 5$, with contractibility function $\rho$ then there is an $\epsilon>0$ such that any $LGC(\rho)$ $n$-manifold within $\epsilon$ of M in Gromov-Hausdorff space is homeomorphic to it.  A natural question is whether we can make this relationship depend solely on $\rho$ and $n$? If the answer were yes, one would obtain a straightforward explanation of the following result of Ferry, \cite{F1}.
\mk\ni
\textbf{Theorem} For every $n$, and contractibility function $\rho$, precompact collections of closed $LGC(\rho)$ Riemannian manifolds in Gromov-Hausdorff space contain only finitely many homeomorphism types.\footnote{The reference cited only proves this in dimensions $n \ne 3$. Recent developments in 3-dimensional topology have corrected this deficiency.}
\mk\ni
We could cover this precompact space by suitable $\epsilon$-balls, refine to get a finite subcover, and each of these balls would contain a unique homeomorphism type.  This strategy is correct with respect to homotopy types, as mentioned above, and for simple homotopy types and rational Pontrjagin classes, as shown in \cite{F1}.
\mk\ni
We will see that, in general, the answer is no.  It can happen that for suitable precompact collections of closed Riemannian manifolds with contractibility function $\rho$, for every $\epsilon$, there are $\epsilon$-balls containing more than one homeomorphism type.
\mk\ni
The most straightforward way to detect this phenomenon is via the symbol of the signature operator on a Riemannian manifold. This lies in $KO_{*}(M)$.  (By work of Sullivan and Teleman, this makes sense for topological manfiolds, except in dimension 4.)  We will show that there exist $\rho$ and arbitrarily close manifolds in a precompact subset of $LGC(\rho)$ whose symbols can differ by an odd torsion element of $KO_{*}(M)$. 
\mk\ni
\textbf{Definition.} We will say that closed manifolds $M$ and $N$ are \textit{deformation equivalent} if there are paths $M_{t}$ and $N_{t}$, $0 \le t <1$ in a precompact subset of Gromov-Hausdorff space consisting of manifolds with contractibility function $\rho$ such that the Gromov-Hausdorff distance between $M_{t}$ and $N_{t}$ goes to zero as $t$ approaches 1. It turns out that this relation is an equivalence relation. A manifold whose class consists of one element is said to be \textit{immutable}.
\mk\ni
\textbf{Theorem 1.} If $M^{m}$, $m \ge 7$, is a closed simply connected manifold such that $\pi_{2}(M)$ vanishes, then there are manifolds which are deformation equivalent  to M in a precompact collection of $LGC(\rho)$-manifolds for some $\rho$ if and only if $KO_{m}(M)$ has odd torsion.  Indeed, for each odd torsion class $\tau$ in $KO_{m}(M)$ there is a unique homotopy equivalence $f:N \to M$ which is realized by a deformation and whose signature operator differs from that of $M$ by $\tau$. This gives many simply-connected examples -- for instance between $S^{3}$-bundles over $S^{4}$.
\mk\ni
	For the  general non-simply connected situation, there are secondary invariants that arise in the problem.  These invariants are related to $\eta$ invariants, except that the familiar Atiyah-Patodi-Singer invariants usually give rise to torsion-free invariants, and the generalization we need must contain torsion information. 
	\mk\ni
	We shall give a complete analysis of the deformation problem for dimensions $\ge 7$ in Theorem 2.8.\footnote{Our construction requires us to embed a certain 3-dimensional compact metric space into $M$, which pushes the dimension to 7.}  We mention here some consequences and examples:
\mk\ni
\textbf{Theorem 2.}  For any $M$, the set of homotopy structures $f:M' \to M$ that are obtainable by deformations in some precompact subset of $LGC(\rho)$ manifolds in Gromov-Hausdorff space defines a subset $S^{CE}(M)$ that is an odd torsion subgroup of the structure group $S(M)$ [40].  It is natural under continous maps of oriented topological manifolds that have dimensions differing by a multiple of 4.
\mk\ni
\textbf{Theorem 3.}  If M has word hyperbolic fundamental group, or has fundamental group that is a lattice is a semisimple Lie group, $S^{CE}(M)$ is finite.
\mk\ni
Theorem 3 depends on the work of Farrell and Jones \cite{FJ}, Bartels and L\"uck \cite{BL}, and  Kammeyer-L\"uck-R\"uping \cite{KLR} on the Borel conjecture. See the discussion following Corollary \ref{fgen}.
\mk\ni
\textbf{Theorem 4.}  There is a closed $M$ such that $S^{CE}(M)$ is infinite.
\mk\ni
	This theorem is in sharp contrast with Ferry's theorem mentioned above.  The resolution of this tension is that for any given $\rho$ only a finite subset (no reason to believe it is a subgroup!) of S(M) occurs.  By varying $\rho$ we obtain this plenitude of deformations.
\mk\ni
	  The fundamental group involved in Theorem 4 is linear, being one of the subgroups of right angled Artin groups, studied in \cite{BB} and \cite{LS}.  The invariant that detects infinitely many homeomorphism types is based on a modification of the theory of higher rho invariants of \cite{We}.  
\mk\ni
The Borel conjecture is currently unresolved in its full generality, so the following corollary to our analysis is especially gratifying.
\mk\ni
\textbf{Theorem 5.} If M is aspherical then $S^{CE}(M) = 0$.
\mk\ni
	We now turn to the second source of motivation, and the area from which the proofs develop.  The theorem of Chapman and Ferry mentioned earlier implies that the limit points which are limits of more than one topological type are not manifolds.  It turns out that they are infinite-dimensional homology manifolds.
\mk\ni
	When these limit points are finite dimensional, they are ANRs, and the work of \cite{Q} implies that there can only be a single topological type in a sufficiently small Gromov-Hausdorff neighborhood.  
\mk\ni
	The possibility of infinite dimensional limit points was established by T. Moore in \cite{M}, based on work of the first author  \cite{Dr} and R. D. Edwards \cite{Wa}.  One key to the construction is the acyclicity of certain Eilenberg-MacLane spaces with respect to mod $p$ $KO$-theory, $p$ odd.  The infinite dimensionality of the cell-like images used in this paper is detected by kernel in the induced map on mod $p$ $KO$-homology.  
\mk\ni
	The close connection of mod $p$ $KO$-theory, for odd primes, to surgery theory is behind our results.  For elements $f:M' \to M$ of $S^{CE}(M)$, we construct infinite dimensional integral homology manifolds $X$ (which will vary by the element), for which $X$ is a cell-like image of both $M'$ and $M$.  This cell-like image of $M$ is constructed by a variant of the method of \cite{Dr}, while the map $M' \to X$ is a consequence of controlled surgery theory.\footnote{Our theorem corrects a mistaken missing finite dimensionality hypothesis in Theorem 3.2.3 in \cite{Q}.}
\mk\ni
	The next section gives information about cell-like maps and classical surgery, stating our main results and setting the stage for the work to follow.  It also contains detailed calculations for several classes of manifolds.  Section 3 contains the details of the construction of suitable cell-like images of manifolds and complexes.  The main theorems are proved in sections 4 and 5, using controlled surgery over the cell-like images.  Finally in section 6 this is related to $LGC(\rho)$ subsets of Gromov-Hausdorff space as well other natural geometric questions (such as the existence of a topological injectivity function for deformations).  The proof of this connection depends on our main theorem that constructs deformations.    At the end, we discuss a modification of the higher $\rho$-invariants that contains enough torsion information to give the examples in Theorem 4.

\mk\ni
Our work leaves open the following question:
\begin{ques}
Can nonhomeomorphic Riemannian manifolds $M$ and $M'$ be deformed to
each other in a precompact subset of Gromov-Hausdorff space, respecting
a contractibility function as above,  while maintaining an upper bound
on volume? Greene and Petersen \cite{GP} have shown that this cannot happen in the presence of an upper bound on volume if the contractibility function is H\"older continuous.\end{ques}

\section{Surgery and cell-like maps}

\ni We begin by giving the formal definition of {\it cell-like map}:

\begin{defn}\label{CE} \

    \begin{enumerate}
        \item [(i)] A compact subset \( X \) of \( \bR^{n} \) is said to
        be \textit{cell-like} if for every open neighborhood \( U \) of \( X \) in
        \( \bR^{n} \), the inclusion \( X \to U \) is nullhomotopic.  This is a
        topological property of \( X \) \cite{La1} and is the \v Cech analogue of
        ``contractible''.  The space $\sin(1/x)$-with-a-bar is an example of a cell-like
        set which is not contractible.
        \item [(ii)] A map $f:Y \to Z$ between compact metric spaces is {\it
        cell-like or CE} if for each $z\in Z$, $f^{-1}(z)$ is cell-like.  The empty
        set is not considered to be cell-like, so cell-like maps must be
        surjective.
    \end{enumerate}
\end{defn}

\noindent  Cell-like maps with domain a compact manifold or
finite polyhedron are weak homotopy equivalences over every open subset of
the range \cite{Ko, La2}.  That is, if \( c:M \to X \) is cell-like, then for every open
\( U \subset X \), \( c|_{c^{-1}(U)}:c^{-1}(U) \to U \) is a weak homotopy
equivalence.  The Vietoris-Begle Theorem implies that the range space of such a cell-like map always has finite cohomological dimension.  If the range has finite covering dimension,
then \( c \) is a homotopy equivalence over every open set.
\mk\ni
LIFTING PROPERTY: Let $f:M\to X$ be a cell-like map with \( M \) an absolute neighborhood retract.\footnote{If $M$ is compact metric with finite covering dimension, $M$ is an absolute neighborhood retract $\equiv$ ANR if and only if it is locally contractible.} Given a space \( W \) with $\dim W<\infty$, $\epsilon>0$, a closed subset $A\subset W$,
a map $g:W\to X$, and a map
$h:A\to M$ with $f\circ h=g|_A$, there is a map $\bar h:W\to M$
extending $h$ such that $g$ is $\epsilon$-homotopic to $f\circ \bar h$ rel $A$:
\[
    \xymatrix{ A \ar[rr]^{h}\ar@{_{(}->}[d]& &M\ar[d]^{f}\\
    W \ar[rr]^{g}\ar[urr]^{\bar h} & & X. }
    \]
See \cite{Ko, La2}
for details.

\begin{defn}
    A homotopy equivalence \( f:N \to M \) is \textit{realized by
    cell-like maps} if there exist a space \( X \) and cell-like maps \(
    c_{1}:N \to X \), \( c_{2}:M \to X \) so that the diagram
    \[
    \xymatrix{ N \ar[rr]^{f}\ar[dr]_{CE}& &M\ar[dl]^{CE}\\
    & X &
    }
    \]
    homotopy commutes. We will also say that $f$ 
    \textit{factors through cell-like maps} and we will call such
    manifolds $N$ and $M$ {\it CE-related}. 
\end{defn}
\mk\ni
In view of the lifting property, every pair of cell-like maps
$c_1:N\to X$, $c_2:M\to X$ induces a homotopy equivalence $f:N\to M$. The induced homotopy equivalence is unique up to homotopy. If $\dim X<\infty$ and $n \ge 4$,  Quinn's
uniqueness of resolutions theorem implies that this homotopy equivalence is
homotopic to a homeomorphism. See \cite{Q}, Prop 3.2.3.
\mk\ni
Two simple homotopy equivalences of manifolds $f_1:N_1\to M$ and
$f_2:N_2\to M$ are called {\it equivalent} if there  is a homeomorphism
$h:N_1\to N_2$ such that $f_{2}\circ h$ is homotopic to $f_{1}$. We
recall that the set $\CS^{s}(M)$ of  equivalence classes of simple homotopy
equivalences $f:N\to M$ is called the set of {\it topological
structures} on $M$. The structure set $\CS^{s}(M)$ is functorial and has
an abelian group structure defined either by Siebenmann periodicity
\cite{KS1} or by algebraic surgery theory \cite{Ra}. Ranicki's
theory gives the induced homomorphism formula for topological
structures \cite{Ra2}:

\begin{prop}\label{Ra}
Let $M^n$ be a closed topological $n$-manifold, $n\ge 5$ and let
$h:M\to N$ be a simple homotopy equivalence, $[h]\in\CS^{s}(N)$.
Then the isomorphism $h_*:\CS^{s}(M)\to\CS^{s}(N)$ is defined by the formula
$$
h_*([f])=[h\circ f]-[h].
$$
\end{prop}
\mk\ni
The structure set $\CS^{h}$ is defined similarly, using homotopy equivalences and 
replacing the relation of homeomorphism by $h$-cobordism. The next proposition shows that 
the homotopy equivalences arising most naturally in this paper are simple. We will omit the 
decorations unless we explicitly wish to study $\CS^{h}$. Similarly, $L$ will be an abbreviation 
for $L^{s}$.

\begin{prop}\label{simple}
A homotopy equivalence $f:M\to N$ that factors through cell-like maps is
a simple homotopy equivalence.
\end{prop}
\begin{proof}
Let $p:M\to X$ and $q:N\to X$ be  cell-like maps such that $f$ is a homotopy
lift of $p$ with respect $q$.
Theorem D of \cite{F3} states that there is  a simple homotopy equivalence
$g:M\to N$ such that $p$ is homotopic to $q\circ g$.
In view of the Lifting Property $f$ is homotopic to $g$. This implies the
equality of the Whitehead torsions:
$
\tau(f)=\tau(g)=0.$ Hence $f$ is a simple homotopy equivalence.
\end{proof}
\mk\ni
By $\CS^{CE}(M)\subset\CS(M)$ we denote the subset of structures realized by
cell-like maps.
\begin{thm}\label{main}
Let $M^n$ be a closed simply connected topological $n$-manifold
with finite $\pi_2(M)$, $n>6$. Then $\CS^{CE}(M)$ is the odd
torsion subgroup of $\CS(M)$.
\end{thm}
\mk\ni
The proof of Theorem 2.5 follows Corollary \ref{cestructure}.
\mk\ni
\begin{brem} On page 531 of \cite{La3}, Lacher asks whether two closed manifolds that admit CE maps to the same space $X$ must be homeomorphic. The theorem above shows that the answer to his question can be ``no'' when $X$ is infinite-dimensional. See Corollary \ref{two} below for an example.\end{brem}
\mk\ni
We recall the
Sullivan-Wall surgery exact sequence \cite{Wall} for closed orientable high-dimensional topological
manifolds:

\begin{equation}\label{SW}
\xymatrix{\cdots \ar[r] & L_{n+1}(\bZ\pi_{1}(M)) \ar@{-->}[r]&\CS(M) \ar[r]^(.4){\eta}& [M,\GTop] \ar[r]^{\theta}&
L_{n}(\bZ\pi_{1}(M)).}
\end{equation}
\ni
The map $\eta$ is called the {\it normal invariant} and
the homomorphism $\theta$ is called the {\it surgery obstruction}. The Sullivan-Wall surgery exact sequence was extended by Quinn and Ranicki\footnote{Our notation differs from Ranicki's in that we've shifted the index on the structure set by one and omitted a bar over $\mathcal S$.}
to the functorial exact sequence of abelian groups below:
\begin{equation}\label{funcSW}
\xymatrix{
\cdots \ar[r] & L_{n+1}(\bZ\pi_{1}(M))\ar[r] &\CS_n(M) \ar[r]^{\eta'}& H_n(M;\bL)\ar[r]^{\theta'} \ar[r] & L_{n}(\bZ\pi_{1}(M))\ar[r] &\cdots}
\end{equation}
where $H_n(M;\bL)=H^0(M;\bL)=[M,\GTop\times\bZ]$, $\CS(M)\subset\CS_n(M)$, and $\eta'|_{\CS(M)}=\eta$.
The homomorphism $\theta'$ is called the {\it assembly map} for $M$.
This sequence is defined and functorial when $M$ is a finite polyhedron. 
This was extended to arbitrary CW complexes in \cite{WW}. See also \cite{We}, page 89.
We write $L_n=L_n(\bZ)$ and recall that
$L_n=\bZ$ if $n=4k$,  $L_n=\bZ_2$ if $n=4k+2$,
and $L_n=0$ for odd $n$. Since $\CS_n(pt)=0$ and in view of the naturality
of the Quinn-Ranicki exact sequence,
the periodic normal invariant $\eta'$ has its image in the kernel $\ker(c_*)$ of
the $\bL$-homology homomorphism  induced by a constant map $c:M\to pt$. Thus,
$\eta'$ factors through
the reduced $\bL$-homology group:
$\bar\eta:\CS_n(M)\to \bar H_n(M;\bL)$.
\mk\ni
In general, Ranicki's algebraic surgery functor gives
us a long exact sequence
\[ \cdots \to \CS_{n}(P,\,Q) \to H_{n}(P,\,Q;\bL) \to
L_{n}(\bZ\pi_{1}P,\,\bZ\pi_{1}Q)\to \cdots \] \vskip .1in
\ni for any CW pair \( (P,\,Q) \).  If \( P\) happens to be a
compact \( n \)-dimensional manifold with nonempty boundary $Q$, then \( \CS_n(P) \) is the usual
rel boundary structure set. \( \CS_n(P) \) differs from the usual geometrically defined structure set by at most a \( \bZ \) in the
closed case.  We also have a long exact sequence \[ \cdots \to
\CS_{n+1}(P,\,Q) \to \CS_{n}(Q) \to \CS_{n}(P) \to \CS_{n}(P,\,Q) \to
\cdots\]

\vskip .1in \ni where for an \( n \)-dimensional manifold with nonempty boundary \(
(P,\,\p P) \), \( \ \CS_{n}(P,\,\p P) \) is  the \textit{not} rel boundary
structure set. 
\mk\ni 
All of these sequences are 4-periodic.
If \( Q \to P \) induces an
isomorphism on \( \pi_{1} \), then \( \CS_{k}(P,Q) \cong H_{k}(
P,\,Q;\bL)\) because the Wall groups \( L_{*}(\bZ\pi_{1}P,\,\bZ\pi_{1}Q )\)
are zero.  Composing this isomorphism
with the boundary map in Ranicki's exact sequence, we
have a homomorphism \(\partial': H_{k+1}( P,\,Q;\bL) \to \CS_{k}(Q)\). For a closed
$n$-manifold there is a split monomorphism 
\begin{equation}\label{SvsSn}
\xymatrix {0\ar[r] &\CS(M)\ar[r]^{i}&\CS_n(M)\ar[r]& \bZ}.
\end{equation}
\mk\ni
To state the main theorem for non-simply connected manifolds we need the following.
\begin{defn}
    If \( K \) is a CW complex, let \( E_2(K) \) be the CW complex obtained from
    \( K \) by attaching cells in dimensions 4 and higher to kill the homotopy
    groups of \( K \) in dimensions 3 and above.  Thus, \( K \subset E_2(K) \),
    \( \pi_{i}(E_2(K))=0 \) for \( i\ge 3 \), and \( E_2(K)-K \) consists of
    cells of dimension \( \ge 4 \). Note that \( E_2(K) \) will not, in general, be a
    finite complex. The space \(E_2(K)\) is called {\it the second stage
    of the Postnikov tower of} \(K\).
\end{defn}

\mk\ni
Let $M$ be a closed $n$-manifold. We denote by
\[ \delta :H_{n+1}(E_2(M),M;\bL) \to \CS(M) \]
the composition:
$$
H_{n+1}(E_2(M),M;\bL) \cong \CS_{n+1}(E_2(M),M)\stackrel{\partial}\to
\CS_n(M)\stackrel{p}\to\CS(M).
$$
where $p$ is any splitting of $i$. The isomorphism above follows from the $\pi - \pi$ theorem.
\mk\ni
Let $\phi:A\to B$ be a homomorphism of abelian groups.
By $\phi^T:\T(A)\to \T(B)$ we denote the restriction
$\phi|_{\T(A)}$ of $\phi$ to the torsion subgroups
and by $\phi_{[q]}:A_{[q]}\to B_{[q]}$ we denote the localization of $\phi$
away from $q$.
\mk\ni
Here is our main theorem for non-simply connected manifolds.
\begin{thm}\label{main2}
Let $M^n$ be a closed topological $n$-manifold, $n>6$.
Then $\CS^{CE}(M)=\im(\delta^T_{[2]})$. In particular,
$\CS^{CE}(M)$ is a subgroup of the odd torsion of $\CS(M)$.\footnote{Using results of \cite{BM} and \cite{AH}, one sees that replacing $E_{2}(M)$ by $E_{k}(M)$, $k \ge 2$, would not not change $\im(\delta^T_{[2]})$.}
\end{thm}
\mk\ni
Since torsion elements of $\CS_{n}(M)$ lie in the kernel of the map 
$\CS_{n}(M) \to \bZ$, $\p$ maps $T(\CS_{n+1}(E_2(M),M))$ into $T(\CS(M)) \equiv T(\CS_{n}(M))$, so 
$\im(\delta^T_{[2]})$ is independent of the choice of $p$. Since the study of $\CS^{CE}(M)$ reduces to an analysis of odd torsion, we can invert 2 in many of our applications. This allows us to omit decorations on $L$-groups and structure sets. We emphasize that this is not a choice on our part. $\CS^{CE}(M)$ is naturally an odd torsion subgroup of $\CS(M)$.
\mk\ni
\begin{cor}\label{fgen} If $L(\pi_{1}M) \otimes \bZ[\frac{1}{2}]$ is finitely generated, then $\CS^{CE}(M)$ is finite.
\end{cor}
\mk\ni
This implies Theorem 3. The Farrell-Jones conjecture for $L(\Gamma) \otimes \bZ[\frac{1}{2}]$ only makes use of the equivariant homotopy theory of $\underline{E\Gamma}= $ the classifying space for proper $\Gamma$-actions. For a lattice, $K\text{\textbackslash} G/\Gamma$ is finite (by the Borel-Serre compactification) and similarly the Rips complex is a suitable space when $\Gamma$ is hyperbolic \cite{MS}. In these cases, the Farrell-Jones conjecture is affirmed (even integrally) in \cite{KLR} and \cite{BL}, respectively.

\begin{cor}\label{induce} Let $f_*:\CS(M)\to\CS(N)$ be the induced homomorphism for
a continuous map $f:M\to N$ between two closed $n$-manifolds, $n>6$.
Then $f_*(\CS^{CE}(M))\subset\CS^{CE}(N)$.
\end{cor}
\begin{proof} We have a commuting diagram

$$
\xymatrix{
\CS_{n+1}(E_{2}(M),M) \ar[r]\ar[d]^{f_{*}}&\CS(M)\ar[d]^{f_{*}}\\
\CS_{n+1}(E_{2}(N),N) \ar[r]&\CS(N)\\
}
$$
from which the result follows immediately.
\end{proof}
\begin{cor} Let $n \ge 6$ and let \( f:N\to M \) be a homotopy equivalence between
closed $n$-manifolds that is realized by cell-like maps. Then $f$ preserves rational Pontrjagin classes.
\end{cor}
\begin{proof} This is Remark 1.7 of \cite{F2}. It follows from the observation that 
$\CS^{CE}(M)$ consists of torsion elements and that the restriction of the normal 
invariant to $\CS^{CE}(M)$ is therefore rationally trivial. Since the rational normal 
invariant measures differences of rational $L$-classes, and therefore rational Pontrjagin classes, 
rational Pontrjagin classes are preserved by structures in $\CS^{CE}(M)$. Crossing with $\bC\bP^{2}$ extends this to all dimensions.
\end{proof}

\begin{cor} Being CE-related is an equivalence relation on closed $n$-manifolds,
$n>5$.
\end{cor}
\mk\ni
{\it Proof.} We prove transitivity. Let $M_1$ be CE-related to
$M_2$ and $M_2$ CE-related to $M_3$. Let $h_1:M_1\to M_2$ and
$h_2:M_2\to M_3$ be corresponding homotopy equivalences. It
suffices to show that the composition $h_2\circ h_1$ is realized by
cell-like maps. In view of Corollary~\ref{induce} we have
$(h_2)_*([h_1])\in\CS^{CE}(M_3)$ and hence by the formula for the
induced homomorphism (Proposition~\ref{Ra} ) we obtain that
$[h_2\circ h_1]=[h_2]+(h_2)_*([h_1])\in\CS^{CE}(M_3)$ .\qed

\mk\ni
We will refer to a manifold that admits a nontrivial deformation as being ``malleable''. 
Manifolds which are not malleable are ``immutable''. In special cases, it is not hard to understand the map \(
H_{n+1}(E_2(M),M;\bL) \to  \mathcal{S}(M)\) well enough to  get
concrete ``immutability'' and ``malleability'' results. We begin with two
typical immutability statements:

\begin{cor}\label{rigid}
    If \( M^{n} \) is a closed manifold with \( n \ge 6  \) and either
    \begin{enumerate}
        \item  [(i)]   \( M \) is aspherical, or

        \item  [(ii)]  \( M \) is homotopy equivalent to a complex
        projective space, or
        \item  [(iii)]  \( M \) is homotopy equivalent to a lens space,
    \end{enumerate}
    \ni  then any
    homotopy equivalence \( f:N \to M \) that factors through cell-like
    maps is homotopic to a homeomorphism.
\end{cor}

\begin{proof}
    If \( M \) is aspherical, then \( M=E_2(M) \) and \(
    H_{n+1}(E_2(M),M;\bL)=0 \), so structures in the image of \(
    H_{n+1}(E_2(M),M;\bL)=0  \) are trivial.
    \mk\ni
    If \( M \) is homotopy equivalent to \(
    \bC P^{k}\), then \( E_2(M)=\bC P^{\infty} \). But \(
    H_{n+1}(\bC P^{\infty},\bC P^{k};\bL)=\lim_{\ell \to \infty}
    H_{n+1}(\bC P^{\ell},\bC P^{k};\bL)\), which has no odd
    torsion, so no nontrivial element of \( \mathcal{S}(M) \) can be the
    image of an odd torsion element. See Lemma \ref{torsion} below.
    \mk\ni
    If $M$ is homotopy equivalent to a $2k-1$-dimensional lens space, 
    then $E_{2}(M)$ is an infinite dimensional lens space, constructed by 
    attaching one cell in each dimension $2k$ and above to $M$. It's 
    straightforward to write down the chain complex for $C_{*}(E_{2}(M),M)$ 
    and compute the integral homology. A quick calculation using the Atiyah-Hirzebruch spectral sequence shows 
    that $H_{2k}(E_{2}(M),M;\,\bL)=\bZ$, so $M$ is immutable.
\end{proof}

\begin{lem}\label{torsion}
    If \( (K,L) \) is a finite CW pair and \( H_{*}(K,L;\bZ) \) has no
    odd torsion, then \( H_{*}(K,L;\bL) \) has no odd torsion.
\end{lem}

\begin{proof}
    For finite CW pairs, \( H_{n}(K,L;\bL) \otimes \mathbb{Q} \cong
    \bigoplus_{k}H_{n-4k}(K,L;\mathbb{Q})\). Comparing this to the
    Atiyah-Hirzebruch spectral sequence gives the result, since there
    can be no nonzero differentials between terms of the form \(
    H_{p}(K,L;L_{4k}) \) on the \( E_{2} \)-page.
\end{proof}

\begin{cor}
It follows that all simply connected manifolds with finite $\pi_{2}$ and no odd torsion in homology are immutable in the sense of Corollary \ref{rigid}.
\end{cor}
\mk\ni
Here is a simple example of malleability.
\mk\ni
\begin{cor}\label{two} There are closed nonhomeomorphic $7$-dimensional
manifolds $M$ and $N$ which are CE-related.
\end{cor}
\begin{proof}
Let $p\ge 7$ be a  prime number. By general position, the Moore complex
$P=S^1\cup_pB^2$ can be PL-embedded in $\bR^5$. Suspending twice 
embeds $P'=S^{3}\cup_{p}B^{4}$ into $\bR^{7}$. Let $Q$ be a regular neighborhood 
of $P'$ in $\bR^{7}$ and let $W = Q \x [0,\,1]$. We set
$\p W = M$, so $M$ is the double of $Q$. $M$ has no ordinary homology
below dimension 3, so $M$ is 2-connected.
By Lefschetz duality, $H_3(W,M)=H^5(W)=H^5(P')=0$ and
$H_4(W,M)=H^4(W)=H^4(P')=\bZ_p$. The exact sequence of the pair
$(W,M)$ turns into:
$$
0\to\bZ_p\to H_3(M)\to \bZ_p\to 0,
$$
which is split by the composition $W = Q \x I \to Q \to M$.
By the Atiyah-Hirzebruch spectral sequence $\bar H_3(M;\bL)\cong \bZ_{p}\oplus \bZ_{p}$.  
Choose a nontrivial $p$-torsion
element $\alpha\in \bar H_7(M;\bL)\cong \bar H_3(M;\bL)$. 
\mk\ni
We will make use of a classical result:

\begin{prop}\label{iso}
For a simply connected closed $n$-manifold  $M$ the reduced normal invariants
$\bar \eta:\CS(M)\to \bar H_n(M,\bL_{\bullet})$ and $\bar \eta':\CS_{n}(M)\to \bar H_n(M,\bL)$ are isomorphisms. Here, $\bL_{\bullet}$ is the 
connective cover of $\bL$. Dualizing, we have an isomorphism $\CS(M) \to [\circM,\,G/TOP]$ and an epimorphism $\CS^{DIFF}(M) \to [\circM,\,G/O]$.
\end{prop}
\begin{proof} Apply the $\pi-\pi$
theorem to  the pair $(M, pt)$ and the pair $(M-\circD, \p D)$, where $D\subset M$ is an $n$-disk.
\end{proof}

\mk\ni
Let
$\beta=\bar\eta'^{-1}(\alpha)$. Since $\beta$ is a torsion element, $\beta \in \CS(M) \subset \CS_{7}(M)$. Thus, by Theorem \ref{main},
$\beta$ defines a homotopy equivalence $f:N\to M$ that belongs to
$\CS^{CE}(M)$. It remains to show that $N$ is not homeomorphic to $M$. 
\mk\ni
We have a diagram:
$$
\xymatrix{
\bZ \ar[r]&\CS^{DIFF}(M) \ar[r]\ar[d] & [\circM;\,G/O]\ar[d]\\
&\CS(M) \ar[r] & [\circM;\,G/TOP] &\ar[l] [\circM;\,BTOP]
}
$$
If we invert 2, and 3, $G$ is 6-connected and $TOP/O$ is 6-connected, so the vertical map on the right gives an isomorphism between $[M;\,G/O]$ and $[M;\,G/TOP]$ which are both equivalent, after the same inversions, to $[M;\,BTOP]$. Since the tangent bundle to $M$ is stably trivial, the topological tangent bundle of $N$ must be stably nontrivial and $N$ cannot be homeomorphic to $M$.
\end{proof}

\mk\ni Thanks to Diarmuid Crowley for pointing out malleable examples in which both $M$ and $N$ are smooth 3-sphere bundles over $S^{4}.$ Diarmuid also points out that $\pi_{n}(PL/O)=0$ for $n \le 6$, \cite{MM}, page 89, so the first Pontrjagin class is a PL invariant. One could use this, together with the Kirby-Siebenmann obstruction to show that our normal invariants are topological invariants in the constructions of this section. The argument above, which is well-known, has the virtue of showing that inverting the primes that appear in stable homotopy groups and those that appear in the orders of groups of smooth homotopy spheres gives topological invariance in any dimension.

\begin{prop}\label{Diarmuid} 
There are nonhomeomorphic $S^{3}$-bundles over $S^{4}$ which are equivalent under deformation.
\end{prop}
\begin{proof}
$\pi_{3}(SO(4))\cong \bZ\oplus\bZ$, so 3-sphere bundles over $S^{4}$ are classified by pairs of integers $(m,n)$ corresponding to elements $m\sigma + n\rho \in \pi_{3}(SO(4))$ with respect to generators $\sigma$, $\rho$ introduced by James and Whitehead and described in \cite{CE}. If $M_{m,n}$ is the sphere bundle corresponding to $m\sigma + n\rho$, we have $H^{4}(M_{m,n})\cong \bZ_{n}$ and the only other nonvanishing cohomology groups are $H^{0}(M_{m,n})\cong H^{7}(M_{m,n})\cong\bZ$. 
\mk\ni
The paper \cite{CE} gives a complete classification of these manifolds up to homotopy equivalence, homeomorphism, and diffeomorphism and includes a computation of normal invariants of homotopy equivalences between nonhomeomorphic manifolds, allowing a complete classification of these manifolds up to deformation. The classification is somewhat lengthy to write down, however, so we content ourselves with discussing homotopy equivalences $f_{p}:M_{0,p}\to M_{12,p}$ in case $p$ is prime, $p\ge 5$.
\mk\ni
According to  Proposition 2.1 of \cite{CE}, these bundles are fiber-homotopy equivalent and the normal invariant of this equivalence is $1 \in \bZ_{p}$. The bundle $M_{0,p}$ is trivial, so the argument of Corollary \ref{two} shows that $f_{p}$ is realized by a deformation, while $M_{0,p}$ and $M_{12,p}$ fail to be homeomorphic.

\mk\ni
We remind the reader that this means that there exist a contractibility
function $\rho:[0,\,R) \to [0,\,\infty)$ and precompact families of
Riemannian metrics $M_{0,p,t},\,M_{12,p,t},\ 0 <t\le 1$ such that
$\rho$ is a contractibility function for each of these metrics and such
that $\lim_{t \to 0}d_{GH}(M_{0,p,t},M_{12,p,t})=0$. Crossing with
$\mathbb{CP}^{2}$, for instance, gives examples with a topological injectivity
function\footnote{See Definition \ref{contr}.}. The full classification of $S^{3}$-bundles over $S^{4}$ up to deformation can be extracted from the computations in \cite{CE}.
\end{proof}

\ni
When $M$ is simply connected, the group $\CS^{CE}(M)$ is finitely generated and in principle it can be computed
by means of Proposition~\ref{iso}. For the example from Corollary \ref{two}, we have:\begin{prop}\label{computation}
$\CS^{CE}(M)\cong\bZ_{p}\oplus\bZ_p.$
\end{prop}
\begin{proof}
$H_{7}(M;\bL_{\bullet})\cong [M;\,\GTop]$ and $\bar H_{7}(M;\bL)\cong [M,*;\,\GTop \x \bZ,*]$, so $$\CS(M)=\CS^{CE}(M)=H_{7}(M;\bL_{\bullet})\cong \bar H_{7}(M;\bL)\cong \bZ_{p}\oplus \bZ_{p}.$$
\end{proof}
\mk\ni
Crossing with a sphere produces further examples of malleability.
\newpage
\begin{cor}\label{cross}\ 

\begin{itemize}
\item [(i.)] If \( f:M' \to M \) is a simple-homotopy equivalence between closed $n$-manifolds with odd order normal
    invariant in \( H_{n}(M;\bL) \), then \( \id_{S^{k}}\x f :S^{k}\x M' \to S^{k}\x M \) factors through cell-like maps, $k \ge 3,\, n+k\ge 7.$
\item [(ii.)] If \( f:M' \to M \) is a homotopy equivalence, not necessarily simple, with odd order normal invariant in \( H_{n}(M;\bL) \) and $k$ is even, then $\id_{S^{k}}\x f$ is $h$-cobordant to a map that factors through cell-like maps, $k \ge 3,\, n+k\ge 7.$
\item [(iii.)] If \( f:M' \to M \) is a homotopy equivalence, not necessarily simple, with odd order normal invariant in \( H_{n}(M;\bL) \) and $k$ is odd, then  \( \id_{S^{k}}\x f :S^{k}\x M' \to S^{k}\x M \) factors through cell-like maps, $k \ge 3,\, n+k\ge 7.$
\end{itemize}
    \end{cor}

\begin{proof}
 (i)   If \( f:M' \to M \) is a simple-homotopy equivalence, then the normal invariant of $[f]$ is a homotopy class of maps $\eta(f): M \to \GTop$. The normal invariants of $\id_{S^{k}}\x f $ and $\id_{B^{k+1}}\x f $ are the composition of projection onto $M$ with $\eta(f)$, so the normal invariants of $\id_{S^{k}}\x f $ and $\id_{B^{k+1}}\x f $ have odd order. Pushing forward to $\bL$ and dualizing, we have a diagram
    {\small\[
    \xymatrix@R.5in@C.1in{
    H_{n+k+1}(B^{k+1} \x M, S^{k}\x M;\bL)\ar[r]^(.52){\cong} &\mathcal S^{s}_{n+k+1}(B^{k+1} \x M, S^{k}\x M)\ar[d]
    \ar[rr]^(.60){\p}&& \mathcal S^{s}_{n+k}(S^{k}\x M) \\
   &\mathcal S^{s}_{n+k+1}(E_2(B^{k+1}\x M),S^{k}\x M)\ar[rru]^{\p}&&\\
    H_{n+k+1}(E_2(S^{k}\x M), S^{k}\x M;\bL)\ar[r]^(.53){\cong}&\mathcal S^{s}_{n+k+1}(E_2(S^{k}\x M),
    S^{k}\x M)\ar[u]_{\cong}\ar@!<350pt,-10pt>^(.6){\p}&}
    \]}
    \ni 
    where both horizontal isomorphisms come from the $\pi-\pi$ theorem and we have used the inclusion-induced homotopy equivalences $E_{2}(M)\cong E_{2}(S^{k}\x M)\cong E_{2}(B^{k+1}\x M)$ when $k\ge 3$. This shows that \( [\id_{S^{k}}\x f]\) is in the image of \( H_{n+k+1}(E_2(S^{k}\x M),
    S^{k}\x M;\bL) \) via the composition of the dashed arrow with the bottom horizontal arrow and that it comes from an odd order element, namely, the image of $\eta(\id_{B^{k+1}}\x f)$ in $\mathcal S_{n+k+1}(E_2(B^{k+1}\x M),S^{k}\x M)$.
\mk\ni
(ii) Consider the diagram above with $\CS^{s}$ replaced by $\CS^{h}$. The homotopy equivalence $f:M' \to M$ satisfies a symmetry $\tau(f)=(-1)^{n-1}\tau(f)^{*}$ which can be seen in the PL case by computing torsions using triangulations and dual triangulations. This shows that $2\tau(f)=(-1)^{n-1}(\tau(f)^{*}+(-1)^{n-1}\tau(f)).$ Torsions of the form $\tau + (-1)^{n-1}\tau$ can be varied away by including into an $h$-cobordism of torsion $\tau$ and retracting to the other end. After crossing with $S^{k}$, $\tau(\id_{S^{k}}\x f)=2\tau(f)$ so $[\id_{S^{k}}\x f]$ therefore lies in the image of $\mathcal S^{s}_{n+k}(S^{k}\x M)$ in $\mathcal S^{h}_{n+k}(S^{k}\x M)$., which is to say that $\id_{S^{k}}\x f$ is $h$-cobordant to a simple homotopy equivalence. The result follows as in case (i).
\mk\ni
(iii) The product formula for Whitehead torsion implies that $\tau(\id_{S^{k}}\x f)=0$ and argument in (i) applies as above.
\end{proof}

\mk\ni
We recall that by  definition a fake lens space of order $p$ is the orbit space 
of a free action of $\bZ_p$ on a sphere. Since simple homotopy equivalent lens spaces are diffeomorphic, the actions giving rise to the fake lens spaces $L'$ below are topologically nonlinear. Explicit constructions of fake lens spaces as quotients of Brieskorn spheres are studied in \cite{Or}.

\begin{cor}\label{lens} There exist a 5-dimensional lens space $L$  and a fake lens space $L'$
such that $L'\x S^3$ and $L\x S^3$ are CE-related 
and $L' \x S^{3}$ and $L \x S^{3}$ are not homeomorphic. 
\end{cor}

\begin{proof} Let $L$ be the lens space $L_{11}(1,1,3)$ in the notation of \cite{Milnor}, p. 403. The first Pontrjagin class of this manifold is zero. We denote by
$\bL_{\bullet}$ the connected cover of the spectrum $\bL$. Thus,
$H_i(pt;\bL_{\bullet})=L_i $ for $i>0$ and $H_i(pt;\bL_{\bullet})=0$ for $i\le 0$.
In the Atiyah-Hirzebruch spectral sequence for the lens
space $L$ the term $E^2_{1,4}=H_{1}(L;H_4(pt;\bL_{\bullet}))=H_1(L)=\bZ_{11}$
survives to $E^{\infty}_{1,4}$ and hence to $H_5(L;\bL_{\bullet})$. Since $L_5(\bZ\bZ_{11})=0$,
this homology class comes from the structure set $\CS(L)$. Thus there is a simple homotopy equivalence
$f:L'\to L$ with nontrivial normal invariant of order $11$. It follows from Proposition~\ref{cross} that $L\x S^3$ and $L'\x S^3$ are CE-related and from Theorem \ref{push} that they deform to each other.
\mk\ni
By Proposition~\ref{cross} the simple homotopy equivalence
$f\x 1_{S^3}:L'\x S^3\to L\x S^3$ also has nontrivial normal invariant of order $11$, so
by the argument of Corollary~\ref{two} the manifolds $L\x S^3$ and
$L'\x S^3$ cannot be homeomorphic, since this first Pontrjagin class is topologically invariant and zero for $L$ and nonzero for $L'$. 
\end{proof}

\begin{rem}
According to \cite{Milnor}, $L_{11}(1,1,4)$ and $L_{11}(1,6,4)$ are homotopy equivalent. The first Pontrjagin class of $L_{11}(1,1,4)$ is zero and the first Pontrjagin class of $L_{11}(1,6,4)$ is three, so they are not homeomorphic and, therefore, not simple-homotopy equivalent. The normal invariant of the homotopy equivalence has odd order, so crossing with $S^{3}$ produces nonhomeomorphic lens spaces crossed spheres that deform to each other.
\mk\ni
In contrast, if $L$ and $L'$ are as in Corollary \ref{lens}, then the induced homotopy equivalence between $L'\x S^{2}$ and $L\x S^{2}$ is not realized by a deformation. $E_{2}(L \x S^{2}) = L^{\infty}\x \CP^{\infty}$, where $L^{\infty}$ is an infinite lens space. The K\"unneth Theorem for ordinary homology together with the Atiyah-Hirzebruch spectral sequence shows that the normal invariant of $f$ does not go to zero in $H_{7}(E_{2}(L \x S^{2}) ;\bL)$, so the structure $[f]$ cannot lie in the image of $H_{8}(E_{2}(L \x S^{2}), \,L \x S^{2};\bL)$.
\end{rem}

\ni
\textbf{Manifolds with finite fundamental group.}
\mk\ni
\begin{prop}\label{finfun}  Let $M$ be a manifold with finite fundamental group such that $\pi_{2}(M)$ is finite.  Then for $n>6$, $$S^{CE}(M^{n}) \cong \Tor(\Ker(KO_{n}(M) \to KO_{n}(K(\pi_{1}M,1))), Z_{(2)}/Z).$$\end{prop}
\mk\ni
\begin{proof} As in the proof of Lemma \ref{E-two},  $E_{2}(M) \to E_{1}(M) = K(\pi_{1}M,1) = B\pi$ induces an isomorphism on $L\wedge M(p)$ homology, so we can use  $K(\pi_{1}M,1)$ and $E_{2}(M)$ interchangeably in our calculations. Also, $\CS(M) \cong \CS(M,*)$.
\mk\ni
To begin, we have a commuting diagram of surgery exact sequences below:

$$\xymatrix@C=1.0pc @M=0pc{
&&&H_{n+1}(E_{2}(M),*;\bL)\ar[d]\\
&&\CS_{n+1}(E_{2}(M),M)\ar[r]^{\cong}\ar[d] & H_{n+1}(E_{2}(M),M; \bL)\ar[d]\\
H_{n+1}(M,*; \bL)\ar[r]\ar[d] & L_{n+1}(\pi_{1}M,e)\ar[r]\ar[d]& \CS_{n}(M,*)\ar[r]\ar[d] & H_{n}(M,*; \bL)\ar[r]\ar[d] & L_{n}(\pi_{1}M,e)\ar[d]\\
H_{n+1}(E_{2}(M),*; \bL) \ar[r] &L_{n+1}(\pi_{1}M,e)\ar[r] & \CS_{n}(E_{2}(M),*)\ar[r]& H_{n}(E_{2}(M),*; \bL)\ar[r]& L_{n}(\pi_{1}M,e)
}$$

\mk\ni
By the $\pi-\pi$ theorem, the top horizontal arrow is an isomorphism. Since we are interested in odd primary behavior, we can invert 2, which replaces $\bL$-homology by $KO$-homology and $E_{2}(M)$ by $B\pi$, where $\pi=\pi_{1}(M)$. This gives us the diagram below at odd primes:
$$\xymatrix{
&&&KO_{n+1}(B\pi,*)\ar[d]\\
&&\CS_{n+1}(B\pi,M)\ar[r]^{\cong}\ar[d] & KO_{n+1}(B\pi,M)\ar[d]\ar[dl]_{\delta}\\
KO_{n+1}(M,*)\ar[r]\ar[d] & L_{n+1}(\pi,e)\ar[r]\ar[d]^{\cong}& \CS_{n}(M,*)\ar[r]\ar[d] & KO_{n}(M,*)\ar[r]\ar[d] & L_{n}(\pi,e)\ar[d]\\
KO_{n+1}(B\pi,*) \ar[r] &L_{n+1}(\pi,e)\ar[r]^{\cong \otimes \bQ} & \CS_{n}(B\pi,*)\ar[r]& KO_{n}(B\pi,*)\ar[r]& L_{n}(\pi,e)
}$$

\mk\ni In \cite{Wall2}, page 2, Wall shows that for $\pi$ finite, $L_{n}(\pi)$ is the direct sum of a free abelian group and a 2-torsion group. By a transfer argument, see \cite{Ad}, the reduced $KO$-homology of $B\pi$ is torsion. An odd torsion element in the kernel of $\CS_{n}(M,*) \to KO_{n}(M,*)$ must be the image of an element $\alpha \in L_{n+1}(\pi_{1}M,e)$ of infinite order. An odd order element of the kernel of $\CS_{n}(M,e) \to KO_{n}(M,*)$ intersect $\im \delta$ would therefore map to an element of infinite order in $\CS_{n}(B\pi,*)$, which is impossible.
\mk\ni
If $\beta \in KO_{n}(M,*)$ is an odd torsion element mapping to zero in $KO_{n}(B\pi,*)$ then $\beta$ is the image of $\beta' \in KO_{n+1}(B\pi,M)$. 
Since $KO_{n+1}(B\pi,*)$ is torsion, $\beta'$ cannot have infinite order and must be odd torsion. This shows that $\CS^{CE}(M,*)$ maps onto the odd torsion in the kernel of 
$KO_{n}(M,*) \to KO_{n}(B\pi,*)$, completing the proof of the proposition.\end{proof}
\mk\ni
\begin{rem} A related result holds for manifolds with abelian fundamental group $\pi = \bZ^{k}\oplus A$. Splitting off infinite cyclic factors using Shaneson's thesis shows that the Wall groups of finitely generated abelian groups are sums of free abelian groups and finite 2-groups. $B\pi = T^{k}\x BA$ and an easy spectral sequence argument shows that the groups $KO_{*}(T^{k},B\pi)$ are torsion, where $B\pi \to T^{k}$ is the projection. The result is a diagram (see below) with the same formal properties as the second diagram in the proof of Proposition \ref{finfun}. The short exact sequence $$ \to \CS_{*+1}(T^{k},M)\to \CS_{*}(M) \to \CS_{*}(T^{k})\to $$ shows that $\CS_{*+1}(T^{k},M)\cong \CS_{*}(M)$. Comparing the long exact $KO$-homology sequences of $(T^{k},M)$ and $(T^{k},B\pi)$ shows that the odd $KO$-homology in the kernel of $KO_{n+1}(T^{k},M)\to KO_{n+1}(T^{k},B\pi)$ is isomorphic to $\CS^{CE}(M).$
$$\xymatrix{
&&&KO_{n+1}(T^{k},B\pi)\ar[d]\\
&&\CS_{n+1}(B\pi,M)\ar[r]^{\cong}\ar[d] & KO_{n+1}(B\pi,M)\ar[d]\ar[dl]_{\delta}\\
KO_{n+2}(T^{k},M)\ar[r]\ar[d] & L_{n+2}(\bZ^{k},\pi)\ar[r]\ar[d]^{\cong}& \CS_{n+1}(T^{k},M)\ar[r]\ar[d] & KO_{n+1}(T^{k},M)\ar[r]\ar[d] & L_{n}(\bZ^{k},\pi)\ar[d]\\
KO_{n+1}(T^{k},B\pi) \ar[r] &L_{n+1}(\bZ^{k},\pi)\ar[r]^{\cong \otimes \bQ} & \CS_{n}(T^{k},B\pi)\ar[r]& KO_{n}(T^{k},B\pi)\ar[r]& L_{n}(\bZ^{k},\pi)
}$$
\end{rem}\qed

\begin{prop} Let $M^{n}$ be a closed manifold, $n \ge 7$, $\pi_{1}(M)=\pi$ and with $\pi_{2}(M)$ zero or finite. If $\pi$ has split injective assembly map away from 2, then $\CS^{CE}(M)$ is isomorphic to the odd torsion subgroup of $H_{n+1}(B\pi,M;\bL)$.
\end{prop}
\begin{proof} Consider the diagram below (away from 2).
$$
\xymatrix{
&&\CS_{n+1}(B\pi,M)\ar[d]\ar[r]^{\cong}  & H_{n+1}(B\pi,M;\bL) \ar[d]\ar[dl]_{\delta}\\
H_{n+1}(M;\bL)\ar[d]\ar[r] & L_{n+1}(\pi)\ar[d]_{id}\ar[r]^{w}  & \CS(M) \ar[d]_{q}\ar[r] & H_{n}(M;\bL)\ar[d]\\
H_{n+1}(B\pi;\bL)\ar[r]_{A}  & L_{n+1}(\pi)\ar[r] \ar@/_/[l]_{j}  & \CS_{n}(B\pi)\ar[r]^{0}\ar@/_/[l]_{s}  & H_{n}(B\pi;\bL)\\
}
$$
The assembly map $A$ is a split monomorphism, so there are splittings $j$ and $s$, as shown. $w \circ s$ is a splitting of $q$, so $\delta$ maps $H_{n+1}(B\pi,M;\bL)$ isomorphically onto a direct summand of $\CS(M)$. The result now follows from Theorem \ref{main2}. A great many torsion-free groups satisfy this version of the Novikov conjecture. See \cite{Hig,STY}.
\end{proof}

\begin{rem} One can unify some of the calculations we have given in this section when the $C^{*}$-algebra assembly map is known to be split injective and one consequently has ``a refined normal invariant'' $\CS(M) \to KO_{n+1}^{\pi}(\underbar E \pi, \widetilde M) [1/2]$ analogous to the map in the preceding proposition given by the  projection $\CS(M)\to H_{n+1}(B\pi,M;\bL)$ (under an $\bL$-theory integral Novikov hypothesis). In that case, $\CS^{CE}(M)$ is isomorphic to the the image of the odd torsion of  $KO_{n+1} (B\pi, M )\cong KO_{n+1}^{\pi} (E\pi, \widetilde M )$ in $KO_{n+1}^{\pi} (\underbar E \pi, \widetilde M) [1/2]$. \end{rem}

\ni	
\textbf{Spherical space forms.}
\mk\ni
We now give a proof of immutability valid for all spherical space forms. We begin by noting that, as in the last section, the map $\CS(M)\to H_{n}(M;\bL)$ is a monomorphism on odd torsion, so $\CS^{CE}(M)\to H_{n}(M;\bL)$ is a monomorphism.

\mk\ni
As above, for any $X$ with free action of a group $G$ with $p$-Sylow subgroup $G_{p}$, the map $X/G_{p} \to X/G$ is split surjective in any $p$-local homology theory, with a splitting induced by the transfer. In particular, the transfer $\tau: KO_{n}(X/G) \to KO_{n}(X/G_{p})$ is split injective on $p$-torsion.
\mk\ni
Now, let $G$ be a finite group acting freely on $S^{n}$ with quotient $M=S^{n}/G$ and let $p$ be an odd prime. The $p$-Sylow subgroup $G_{p}$ of $G$ must be cyclic, so, as Wall observes, $L=S^{n}/G_{p}$ has the homotopy type of a linear lens space. 
By Wall, \cite{Wall}, Chapter 14, the structure group of an odd lens space is torsion free. See also \cite{We}, pp 110-111. This implies that $S^{CE}(L)$ is trivial. The diagram below then shows that the $p$-torsion in $\CS^{CE}(M)$ must be trivial. The transfer on the left arises from the observation that a map covering a CE map is also CE.

$$
\xymatrix{ 
\CS_n^{CE}(L) \ar[r] & \CS_{n}(L)\ar[r] & H_{n}(L;\bL)\\
\CS_n^{CE}(M) \ar[r]^{1-1}\ar[u]_{\tau}\ar @/_1.0pc/[rr]_{1-1} & \CS_{n}(M)\ar[r]\ar[u]_{\tau} & H_{n}(M;\bL)\ar[u]_{\tau}^{1-1}
}$$
\mk\ni
Repeating for each odd $p$, it follows that $S^{CE}(M)$ is trivial.
\begin{prop}[Proof of Theorem 4] $\CS^{CE}(M)$ can be infinite.\end{prop}

\begin{proof} Let $M(\bZ/p,\,n)$ be a Moore space, $p$ an odd prime. Triangulate $M(\bZ/p,\,n)$ as a flag complex $L$ and form the Bestvina-Brady group $\pi=H_{L}$ as in \cite{BB}, \cite{LS}. For $n \ge 3$, $B\pi$ has free, finitely generated homology through dimension $n$ and its homology in dimension $n+1$ contains an infinite sum of $\bZ/p$'s. See Corollaries 8 and 9 in \cite{LS}. It also follows that $H_{n+2}(\pi)$ is finitely generated free abelian and all higher homology is zero, whence it follows immediately from the Atiyah-Hirzebruch spectral sequence that $H_{n}(\pi,\bL)$ contains an infinite sum of $\bZ/p$'s. By periodicity, the same is true for $H_{n-4k}(\pi,\bL)$ for any $k$.
\mk\ni
Let $K$ be the 3-skeleton of $B\pi$, which is finite for $n \ge 3$. Embed $K$ in $\bR^{m+1}$, $m\ge 8$, and let $M^{m}$ be the boundary of a regular neighborhood. Consider the diagram
$$
\xymatrix{
&&\CS_{m+1}(B\pi,M)\ar[d]_{\cong}\ar[r]^{\cong}  & H_{m+1}(B\pi,M;\bL)\ar[dl]_{\delta}\\
&& H_{m}(M;\bL)\\
}
$$

\ni By the Borel Conjecture for Bestvina-Brady groups \cite{BL}, the vertical arrow on the left is an isomorphism, so $\delta$ is an isomorphism. Since $\pi_{2}(M)=0$, we have $B\pi=E_{2}(M)$, so Theorem \ref{main2} tells us that $\CS^{CE}(M)$ is the image of the odd torsion under the map $\delta$. Let $m+1 = n-4k$, $k \ge 1$. By our construction, $H_{m+1}(B\pi;\bL)$ contains infinitely generated $p$ torsion. Since $H_{m+1}(M;\bL)$ is finitely generated, $H_{m+1}(B\pi,M;\bL)$ contains infinitely generated $p$ torsion and $\CS(M)$ contains infinitely generated $p$ torsion in the image of $\delta$.
\end{proof}

\ni Each structure $[\alpha]$ above is represented by a homotopy equivalence $\alpha:M_{\alpha}\to M$.

\begin{prop} There are infinitely many nonhomeomorphic manifolds $M_{\alpha} \in \CS^{CE}(M)$ when $m \ge 8$.
\end{prop}
\begin{proof}
If $m\ge 8$, the manifold $M$ has a handle decomposition with no handles in the middle dimension. As in \cite{We2}, this allows us to define an absolute ``higher $\rho$ invariant'' for our manifolds $M_{\alpha}$ in a quotient group of $L_{m+1}(\pi)$. Since the construction in \cite{We2} was rational and we are interested in torsion phenomena, we will review the construction. 
\mk\ni
Let $N^{m}$ be a closed, oriented $m$-manifold such that the $\bZ\pi_{1}N$-chain complex $C_{*}(N)$ is chain-homotopy equivalent to chain complex of finitely generated projective $\bZ\pi_{1}N$-modules $\{P_{i}\}$ with $P_{i}=0$, $i=[m/2]$. Following Hausmann, we call such manifolds \textit{anti-simple}. Let $P^{<i}_{*}$ be the truncation of $P_{*}$. There is a chain retraction $P_{*}\to P^{<i}_{*}$ and $(P^{<i}_{*},P_{*})$ is a symmetric algebraic Poincar\'e pair. This is well-defined in that if $Q_{*}$ is a chain complex of finitely generated projective modules chain-homotopy equivalent to $P_{*}$ with $Q_{i}=0$, then there is a chain-homotopy equivalence of pairs $(P^{<i}_{*},P_{*}) \sim (Q^{<i}_{*},Q_{*})$. If the manifold $N$ is the boundary of an oriented manifold $W^{m+1}$ with a map to $B\pi_{1}N$ extending $N \to B\pi_{1}N$, the pair $(W,N)$ gives us another symmetric algebraic Poincar\'e pair over $\bZ\pi_{1}N$ and we can paste the two together along $P_{*}$ to get a closed $(m+1)$-dimensional symmetric algebraic chain complex and, therefore, an element of $L^{m+1}(\pi_{1}(N))$. Two such coboundaries of $N$ define an element $\omega$ of $\Omega_{m+1}(B\bZ\pi_{1}N)$, so our element of $L^{m+1}(\pi_{1}(N))$ is well-defined up to the image of the map $\Omega_{m+1}(B\bZ\pi_{1}N) \to L^{m+1}(\pi_{1}(N))$ that sends each element to its symmetric signature. The resulting element in $L^{m+1}(\pi_{1}(N))/\Omega_{m+1}(B\bZ\pi_{1}N)$ is the higher $\rho$ invariant of $N$.
\begin{rem}
One can define the higher $\rho$ invariant without assuming explicitly that $M^n$ bounds if one inverts the torsion present in Witt bordism in dimension $n$ of $B\pi$.  (This follows from the argument in \cite{We}.)  Thus, for the Bestvina-Brady groups used here, in low dimensions one need not invert anything and one has a more general higher $\rho$ invariant available to distinguish homotopy equivalent anti-simple manifolds with Bestvina-Brady fundamental groups.
\end{rem}
\ni Returning to our manifold $M$, let $P_{*}$ be the $\bZ\pi$-chain complex obtained by gluing together two copies of a handle decomposition of a regular neighborhood of $K$ in $\bR^{m}$. We have $P_{i}=0$ for $3<i<m-3$. Since the manifolds $M_{\alpha}$ are homotopy equivalent to $M$, they are also anti-simple. Since they are obtained from $M$ by Wall realization, they are cobordant to $M$, so they bound and have higher $\rho$ invariants. By construction, the higher $\rho$ invariant of $M_{\alpha}$ differs from the higher $\rho$ invariant of $M$ by the image of $\alpha''$ in $L^{m+1}(\pi)/\Omega_{m+1}(B\pi)$. Since $B\pi$ has finite $m+1$-skeleton, the Atiyah-Hirzebruch spectral sequence shows that $\Omega_{m+1}(B\pi)$ is finitely generated. Since the collection of $\alpha$'s in $L^{m+1}(\pi)$ is infinitely generated, there are infinitely many nonhomeomorphic $M_{\alpha}$'s.
\end{proof}

\mk\ni Next we show that, by itself, infinite odd torsion in the $L$-group does not suffice to produce infinitely many deformable manifolds.
\mk\ni
\begin{prop}There is a closed manifold $M^{n}$ such that $L_{n+1}(\pi)$ has infinitely generated odd torsion but $\CS^{CE}(M)=0$.
\end{prop}
\begin{proof}
Let $A$ be the universal finitely presented acyclic group of \cite{BDM} and let $\Sigma$ be a homology sphere with fundamental group $A$. Since $H_{1}(A)=H_{2}(A)=0$, such a homology sphere exists in dimensions $\ge 5$ by a well-known theorem of Kervaire \cite{K}. The surgery exact sequence for $\Sigma$ is
$$
\xymatrix{
H_{n+1}(\Sigma;\bL)\ar[r]&L_{n+1}(A)\ar[r]&\CS(\Sigma)\ar[r]&H_{n}(\Sigma;\bL)
}
$$
Inspection of this sequence gives us $$\tilde L_{n+1}(A)\cong \CS(\Sigma),$$ where $\tilde L_{n+1}(A) =  L_{n+1}(A)/L_{n+1}(e)$.
Now consider the commutative diagram of topological surgery exact sequences
$$
\xymatrix{
0\ar[r] &\tilde L_{n+1}(A) \ar[r]^{\cong}\ar[d]_{\cong} & \CS(\Sigma)\ar[r]\ar[d] &0\\
0=\bar H_{n+1}(BA;\bL)\ar[r]&\tilde L_{n+1}(A)\ar[r]&\CS(BA)&
}
$$
If an element of $\CS(\Sigma)$ goes to 0 in $\CS(BA)$, then it comes from an element of $\tilde L_{n+1}(A)$ which maps to a nonzero element of $\CS(BA)$, yielding a contradiction. Since elements of $\CS^{CE}(\Sigma)$ must die in $\CS(BA)$, $\CS^{CE}(\Sigma)=0$. If $L_{n+1}(A)$ contains infinitely generated odd torsion, we are done. Otherwise, let $\pi$ be a finitely presented group such that $L_{n+1}(\pi)$ has infinitely generated odd torsion and consider the amalgamated free product $\Gamma = A *_{\pi}(\pi \x \bZ/2)$. $\Gamma$ is $\bZ[1/2]$-acyclic. Using the isomorphism (away from 2) 
$$L(\pi \x \bZ/2)\cong L(\pi)\x L(\pi)$$
and Cappell's Mayer-Vietoris sequence \cite{Ca1}, we have
$$
\xymatrix{
0\ar[r] & \tilde L_{n}(A) \oplus \tilde L_{n}(\pi)\ar[r] & \tilde L_{n}(\Gamma)\ar[r] & 0\\ 
}
$$
\ni after inverting 2. This shows that $\Gamma$ has infinite odd torsion in $L$-theory with $H_{1}(\Gamma)=\bZ/2.$ Suspending once to $A*_{\Gamma}A$ kills $H_{2}$, so a $\bZ[1/2]$ version of Kervaire's theorem produces a $\bZ[1/2]$-homology sphere with fundamental group $A*_{\Gamma}A$. Cappells's theorem mod finitely generated odd torsion shows that the $L$-theory of $L_{n-1}(A*_{\Gamma}A)$ contains infinite odd torsion, and we can complete the argument as above.
\end{proof}

\section{ Cell-like maps that kill $\bL$-classes}
\mk\ni
We use the notation  $KO_*(X)=H_*(X;KO)$ for periodic $KO$-homology and we use $\KO$ to stand  for reduced $KO$ homology.
We need the following facts \cite{AH}, \cite{BM}, \cite{Y}.
\begin{thm}\label{AH}  If $p>1$ is an integer and $n \ge 3$, $\KO(K(\pi,n);\bZ_p)=0$
for any group $\pi$. If $\pi$ is finite,  $\KO(K(\pi,n))=0$, for $n\ge 2$.
\end{thm}
\mk\ni
Let $M(p)$ denote the $\bZ_p$ Moore spectrum. For odd $p$, we have a chain of homotopy equivalences of spectra
$$
\KO\wedge M(p)\sim \KO[\frac{1}{2}]\wedge M(p)\sim\bL[\frac{1}{2}]\wedge M(p)\sim
\bL\wedge M(p).
$$
This implies the following:
\begin{cor}\label{for L} Let $p$ be odd, then
$\bar H_*(K(\bZ_p,2);\bL\wedge M(p))=0$ where $\bL\wedge M(p)$
is $\bL$-theory with coefficients in $\bZ_p$.
\end{cor}

\mk\ni
We recall that for an extraordinary homology theory given by a spectrum $\bE$
of CW complexes there are Universal Coefficient Formulas
for coefficients $\bZ_p$ and $\bQ$:

\begin{equation}\label{UCT}
 0\to H_n(K;\bE)\otimes\bZ_p\to H_n(K;\bE\wedge M(p))\to
\Tor(H_{n-1}(K;\bE),\,\bZ_p)\to 0
\end{equation}
and
$$
H_n(K;\bE_{(0)})=H_n(K;\bE)\otimes\bQ.
$$
Here $\Tor(H,\,\bZ_p)=\{c\in H\mid pc=0\}$ and $\bE_{(0)}$ denotes the localization
at 0.
Let $X=\varprojlim\{K_i\}$ be a compact metric space presented as the
inverse limit of finite polyhedra. By
$\check H_*(X;\bE)=\varprojlim\{H_*(K_i,\bE)\}$ we denote the
\v Cech $\bE$-homology. The Steenrod homology $H_n^{st}(X;\bE)$ of $X$ fits into the following
exact sequence
$$
0\to{\lim}^1\{H_{n+1}(K_i;\bE)\}\to H_n^{st}(X;\bE)\to\check H_n(X;\bE)\to 0.
$$
If $H_k(pt; \bE)$ is finitely generated for each $k$, the Mittag-Leffler condition holds with rational or finite coefficients, so we have
$$
H_n^{st}(X;\bE\wedge M(p))=\check H_n(X;\bE\wedge M(p))\ \ \text{and}\ \ \
H_n^{st}(X;\bE_{(0)})=\check H_n(X;\bE_{(0)}).
$$

\mk\ni
In the case of $\bZ_p$-coefficients we obtain an exact sequence
which is natural in $X$:

\begin{equation}\label{Zp}
 0\to\varprojlim(H_n(K_i;\bE)\otimes\bZ_p)\to H_n^{st}(X;\bE\wedge M(p))
\stackrel{\phi'}\to\Tor(\check H_{n-1}(X;\bE),\,\bZ_p)\to 0.
\end{equation}

\begin{lem}\label{E-two}
Let $M$ be a simply connected finite complex with finite $\pi_2(M)$. Then
for every element $\gamma\in H_k(M;\bL)$ of odd order $p$
there exists an odd torsion element $\alpha\in H_{k+1}(E_2(M),M;\bL)$ such that
$\partial(\alpha)=\gamma$ where $\p$ is the connecting homomorphism
in the exact sequence of the pair $(E_2(M),M)$.
\end{lem}
\begin{proof} Note that $E_2(M)=K(\pi_2(M),2)$.
\mk\ni
If $\pi_2(M)=0$, the space $E_2(M)$ is contractible and the lemma is
trivial.
\mk\ni
If $\pi_2(M)$ is torsion, then by Corollary~\ref{for L},
$\bar H_*(E_2(M);\bL\wedge M(p))=0$.
Then by the Universal
Coefficient diagram

\small{
\begin{equation}\label{UCTbig}
{\xymatrix@C.2in{
H_{k+2}(E_2(M),M;\bL\wedge M(p))\ar[r]^{epi}\ar[d]^{\p}_{iso} &
\Tor(H_{k+1}(E_2(M),M;\bL),\,\bZ_p)\ar[r]^(.58){mono} &
H_{k+1}(E_2(M),M;\bL)\ar[d]^{\p}\\
H_{k+1}(M;\bL\wedge M(p))\ar[r]^{epi} &  \Tor(H_k(M;\bL),\,\bZ_p)\ar[r]^{mono} &
H_k(M;\bL)\\
}}
\end{equation}
}
\normalsize{we obtain the required result.}
\end{proof}
\mk\ni
The following proposition is proven in \cite{Wa} Appendix B.
\begin{prop}\label{Wa}
Let $E$ be a CW complex with trivial homotopy groups $\pi_i(E)=0$,
$i\ge k$ for some $k$, and let $q:X\to Y$ be a cell-like map
between compacta. Then $q$ induces a bijection of the homotopy
classes $q^*:[Y,E]\to[X,E]$.
\end{prop}
\mk\ni
Let $q:M\to X$ be a cell-like map. According to Proposition
\ref{Wa}, for every map $h:M\to E_2(M)$ there is a map $g:X\to
E_2(M)$ such that $g\circ q$ is homotopic to $h$. In particular,
there is an induced map $\tilde g:M_q\to M_h$ between their mapping
cylinders, $\tilde g|_M=id_M$, $\tilde g|_X=g$. We apply this when
$h$ is the inclusion $j:M\subset E_2(M)$ and denote the induced
map by $i:M_q\to M_j$. Denote by
$$i_*:H_*^{st}(M_q,M;\bL)\to H_*(E_2(M),M;\bL)$$ the induced homomorphism
for Steenrod \( \bL \)-homology groups \cite{F1}, \cite{KS}.

\begin{prop}\label{existce}
    Let \( M^{n} \) be a closed connected topological \( n \)-manifold,
    \( n \ge 6
    \), let $p$ be odd, and let \( \beta \in H_{*}(E_2(M), M; \bL\wedge M(p)) \),
     then there exist a cell-like map \( q:M \to X \) and an
    element \( \widehat \beta \in H_{*}^{st}(M_q, M; \bL\wedge M(p) )\)
    such that \( i_{*}(\widehat \beta)=\beta \).
\end{prop}
\mk\ni
The map $q:M \to X$ is a weak homotopy equivalence but is not a homotopy equivalence. 
The proof of Proposition \ref{existce} will follow Lemma \ref{decomposition}.

\begin{prop}\label{existceappl}
    Let \( M^{n} \) be a closed connected topological \( n \)-manifold,
    \( n \ge 6
    \). If \( \alpha \in H_{*}(E_2(M), M; \bL ) \) is an odd torsion
    element, then there exist a cell-like map \( q:M \to X \) and an
    odd order
    element \( \widehat \alpha \in H_{*}^{st}(M_q, M; \bL )\) such that \(
    i_{*}(\widehat \alpha)=\alpha \).
\end{prop}

\begin{proof}
Let $\alpha\in H_k(E_2(M),M;\bL)$ be an
element of order $p$ where $p$ is odd.
Then by the universal coefficient formula, there is an epimorphism
$$\phi: H_{k+1}(E_2(M),M;\bL\wedge M(p))\to \Tor(H_k(E_2(M),M;\bL),\,\bZ_p).$$
Note that $\Tor(H,\,\bZ_p)=\{c\in H\mid pc=0\}$ so that
there is an  inclusion $\Tor(H,\,\bZ_p)\subset H$ which is natural in $H$.
Thus, $\alpha\in
\Tor(H_k(E_2(M),M;\bL),\,\bZ_p)$.
Hence, there is
an element $\beta\in H_{k+1}(E_2(M),M;\bL\wedge M(p))$ such that
$\phi(\beta)=\alpha$. By Proposition \ref{existce} there exist a cell-like map
$q:M\to X$ and an element $\widehat\beta$ such that
\( i_{*}(\widehat \beta)=\beta \). The commuting diagram of universal coefficient formulas gives (see Equation (\ref{UCT})):
$$
\xymatrix{H_{k+1}^{st}(M_q,M;\bL\wedge M(p))\ar[r]^{\phi'}\ar[d]^{i_*}&  \Tor(H_k^{st}(M_q,M;\bL),\,\bZ_p)
\ar[r]^{\subset}\ar[d] & H_k^{st}(M_q,M;\bL)\ar[d]^{i_*}\\
H_{k+1}(E_2(M),M;\bL\wedge M(p))\ar[r]^{\phi} &  \Tor(H_k(E_2(M),M;\bL),\,\bZ_p)\ar[r]^{\subset} &
H_k(E_2(M),M;\bL)\\
}
$$
which implies that $i_*(\hat\alpha)=\alpha$ where $\widehat\alpha=\phi'(\widehat\beta)$
is an element of order $p$.
\end{proof}
\mk\ni
\begin{rem}
 By Proposition \ref{Wa} a cell-like map induces a rational isomorphism on $\bL$-homology.
Therefore, $H_*^{st}(M_q,M;\bL)$ is a torsion group.
\end{rem}

\begin{thm}\label{1-connected}
Let $M^n$ be a closed simply connected topological \( n \)-manifold, $n\ge 6$,
with $\pi_2(M)$ finite.
Then for every odd torsion element $\gamma\in H_*(M;\bL)$
there exist $X$ and a cell-like map $q:M\to X$ such that $q_*(\gamma)=0$.
\end{thm}
\begin{proof}
By Lemma \ref{E-two} there is an odd torsion element
$\alpha\in H_{*}(E_2(M), M; \bL)$
such that $\p(\alpha)=\gamma$.
By Proposition \ref{existceappl} there exists a cell-like map $q:M\to X$ and
an element $\widehat\alpha\in H_{*}^{st}(M_q, M; \bL )$ such that
$i_*(\widehat\alpha)=\alpha$. Then the
commutative diagram
\[
\xymatrix{
H_{*+1}^{st}(M_q,M;\bL) \ar[r]\ar[d]^{i_*} & H_{*}(M;\bL) \ar[r]^{q_*}\ar[d]^{=}
& H_*(X;\bL)\ar[d]\\
H_{*+1}(E_2(M),M;\bL) \ar[r] & H_{*}(M;\bL) \ar[r] &
H_*(E_2(M);\bL)\\
}
\]
implies that $q_*(\gamma)=0$.
\end{proof}
\begin{rem}
 Without the finiteness assumption on  $\pi_2(M)$ one can
show that $q$ kills
an element $\gamma\otimes 1_{\bZ_p}$ with $\bZ_p$ coefficients.
\end{rem}

\mk\ni
We recall that the {\it cohomological dimension} of a topological space $X$
with respect to the coefficient group $G$ is the following number:
$$
\dim_GX=\max\{n\,|\, \check H^n(X,A;G)\ne 0\ \text{for some closed}\ 
A\subset X\}.
$$
\mk\ni
A map of pairs $f:(X,L)\to (Y,L)$ is called {\it strict} if
$f(X-L)=Y-L$ and $f|_L=id_L$.
\mk\ni
The following theorem
is taken from \cite{Dr2} (Theorem 7.2). For $G=\bZ$ it
can be found in \cite{DFW}.

\begin{thm}\label{example}
    Let $\tilde h_*$ be a reduced generalized homology theory.
    Suppose that $\tilde h_*(K(G,n))=0$ for some countable
    abelian group $G$. Then for any finite polyhedral pair
    $(K,L)$ and any element $\alpha\in\tilde h_*(K,L)$ there is
    a compactum $Y\supset L$ and a strict map $f:(Y,L)\to (K,L)$ such that
    \begin{enumerate}
        \item[(i)] $\dim_{G}(Y-L)\le n$;
            \item[(ii)] $\alpha\in \im(f_*)$.
    \end{enumerate}
\end{thm}
\mk\ni
The following is a relative version of  Theorem 3 from \cite{Dr1}.

\begin{lem}\label{resolution}
    Let  $(Y,L)$  be a pair of compacta such that
$$
\dim_{\bZ_p}(Y-L)\le 2  \ \ \text{and} \ \ \ \dim_{\bZ[\frac{1}{p}]}(Y-L)\le 2.
$$
Then there is a strict cell-like map $g:(Z,L)\to (Y,L)$ such that  $\dim(Z-L)\le 3$ and
     $\dim(Z-L)^2\le 5$.
\end{lem}

\begin{lem}\label{retraction}
Let $(Z,M)$ be a compact pair such that $\dim(Z-M)^2\le 2n-1$
and let $M$ be a manifold of dimension $2n$.
Suppose there is a retraction $\rho:Z\to M$.
Then $\rho$ is homotopic $\rel M$ to a retraction
$r:Z\to M$ with $r|_{(Z-M)}$ one-to-one.
\end{lem}
\begin{proof}
The condition $\dim X^2\le 2n-1$ for a compact metric space $X$
implies that every continuous map $\phi:X\to M$ to a $2n$-dimensional
manifold can be approximated by an embedding \cite{DRS}, \cite{Sp}.
Moreover, the space of embeddings $\Emb(X,M)$ is a dense $G_{\delta}$
in the space of mappings $\textrm{Map}(X,M)$. The same argument shows that
under the condition $\dim(Z-M)^2\le 2n-1$
the space of retraction-embeddings $$\textrm{RetEmb}(Z,M)=\{f:Z\to M\mid f|_M=id_M, f|_{Z-M}\ \text{is one-to one}\}$$
is dense in the space
of all retractions $\textrm{Ret}(Z,M)$.
\end{proof}
\mk\ni
The following lemma is proven in \cite{DFW} Lemma 3.7.
\begin{lem}\label{decomposition}
Let $Z$ be a compact and
$r:Z\to M$ be a retraction with $r|_{(Z-M)}$ one-to-one.
Let $g:(Z,M)\to (Y,M)$ be a continuous map which is the identity over $M$.
Then the decomposition of $M$ whose  nondegenerate elements are
$r(g^{-1}(y))$ is upper semicontinuous.
\end{lem}\qed

\ni
{\it Proof of Proposition~\ref{existce}.} We consider the generalized homology theory $h_*=\bL\wedge M(p)$,
i.e., $\bL$-theory with coefficients in $\bZ_p$. Let $\beta\in h_{k+1}(E_2(M),M)$.
There is a finite complex $K$, $M \subset K \subset E_2(M)$, and an element $\gamma\in h_*(K,M)$
such that $\gamma$ is taken to $\beta$ by the inclusion homomorphism.
\mk\ni
Note that $\tilde h_*(K(\bZ[\frac{1}{p}],2))=0$ since the $\bL$-theory of this
space  is $p$-divisible. Taking into account Corollary \ref{for L} and
Theorem \ref{AH} we can state that $\tilde h_*(K(G,2))=0$ for
$G=\bZ_p\oplus\bZ[\frac{1}{p}]$. Then we apply Theorem \ref{example}
to $(K,M)$ and $\gamma$ with this $G$ to obtain $f:(Y,M)\to (K,M)$
satisfying the conditions (i)-(ii) of Theorem \ref{example}.
Condition (i) allows us to apply Lemma \ref{resolution}
to obtain a cell-like map $g:(Z,M)\to (Y,M)$ with $\dim(Z-M)\le 3$ and $\dim(Z-M)^2\le 5$.
\mk\ni
Because $E_2(M)-M$ has no cells of dimension $\le 3$, there is a homotopy of $f\circ g$ rel $M$ that sweeps $Z-M$ to
$M$. Thus, $f\circ g$ is homotopic to a retraction $\rho:Z\to M$.
By Lemma \ref{retraction}, $f\circ g$ is homotopic rel $M$ to a
retraction $r:Z \to M$ which is one-to one on $Z-M$. By Lemma
\ref{decomposition} the decomposition of $M$ into $r(g^{-1}(y))$
and singletons defines a cell-like map $q:M\to X$ such that there
is a commutative diagram
$$
\xymatrix{
Z \ar[rr]^r \ar[d]^g && M \ar[d]^q\\
Y \ar[rr]^{r'} && X\\
}
$$
\mk\ni
By Proposition \ref{Wa} there is a map $g':X\to E_2(M)$ such that
$g'\circ q$ is homotopic to the inclusion $M\subset E_2(M)$.
Hence $f\circ g\sim r\sim g'\circ q\circ r=g'\circ r'\circ g$.
Since $g$ is cell-like,  the map $f$ is homotopic to
$g'\circ r'$ by Proposition \ref{Wa}. Then there is a homotopy commutative diagram of the mapping
cylinders
$$
\xymatrix{
M_{j'} \ar[rr]^{r'}\ar[d]^f && M_q \ar[dll]^{i}\\
M_j\\
}
$$
where $j:M\to E_2(M)$ and $j':M\to Y$ are the embeddings.
For Steenrod $h_*$-homology this gives us the following diagram:
$$
\xymatrix{\tilde h_*(Y,M) \ar[r]^{=}\ar[d]^{f_*}  &\tilde h_*(M_{j'},M)\ar[r]^{r'_*}
\ar[d] & \tilde h_*(M_q,M) \ar[dl]^{i_*}\\
\tilde h_*(E_2(M),M) \ar[r] &\tilde h_*(M_j,M)\\
}
$$
By condition (ii) of Theorem \ref{example} there is
$\gamma'\in\tilde h_*(Y,M)$ such that $f_*(\gamma')=\gamma$.
Then $i_*(\widehat\beta)=\beta$
where $\widehat\beta=r'_*(\gamma')$.
\qed

\

\section{Continuously controlled topology and cell-like maps of simply connected
manifolds}
\mk\ni
We recall that a map of pairs $f:(Z,Y)\to (Z',Y)$ is strict if
$(Z-Y)\subset Z'-Y$ and $f|_Y=id_Y$.
 A proper homotopy $f_t:Z\to X$ which is strict at each level is called {\it strict} if the homotopy
$\bar f_t:(\bar Z,Y)\to(\bar X,Y)$ is continuous.
\mk\ni
Let $X$ be a locally compact space compactified by a compact {\it corona}
$Y=\bar X-X$.
A proper map $f:Z\to X$
is a {\it strict homotopy equivalence} if there is a proper map
$g:X\to Z$ such that
$ g\circ f$ and $ f\circ g$ are strictly
homotopic to $id_{\bar Z}$ and
$id_{\bar X}$ respectively where $Z$ is given a compactification as above.

\begin{defn} \
    \begin{enumerate}
        \item [(i)] Let $X$ be an open manifold and let
        $Y$ be the compact corona of a compactification
        $\bar X$ of  $X$. Two strict homotopy
        equivalences $f:W\to X$
        and $f':W'\to X$ are
        {\it equivalent} if there is a  homeomorphism
        $h:W\to W'$ such that $f=f'\circ h$.
        \item [(ii)]  The set of the equivalence classes
        of strict homotopy equivalences of manifolds is called
        the set of {\it continuously controlled structures}
        on $X$ at $Y$ and it is denoted by
        $\CS^{cc}(\bar X,Y)$.
    \end{enumerate}
\end{defn}
\mk\ni
We note that if $\tilde X$ is another compactification  of $X$
with compact corona $Y'$ such that there is a continuous strict map $\phi:\bar X\to\tilde X$ which is
the identity on $X$, then there is a map $\phi_*:\CS^{cc}(\bar X,Y)
\to \CS^{cc}(\tilde X,Y')$.

\begin{defn} A pair $(X,Y)$ is said to be \textit{locally 1-connected at $Y$} if for each $y \in Y$ and neighborhood $U$ of $y$ in $X$ there is a smaller neighborhood $V$ of $y$ in $X$ so that the inclusion-induced map $\pi_1(V-Y) \to \pi_1(U-Y)$ is zero.\end{defn}

\begin{prop}\label{ccsurgery}
Let $X$ be a simply connected open  manifold of dimension $n\ge 5$
compactified by a compact corona $Y$ in such a way that  the pair $(\bar X,Y)$
is locally 1-connected.
Then there is a surgery exact sequence
$$
\dots\to \bar H_{n}(Y;\bL)\to\CS^{cc}(\bar X,Y)\to [X,\GTop]\to \bar
H_{n-1}(Y;\bL)
$$
which is natural with respect to maps between coronas, as above.
Here $\bar H_*(-;\bL)$ is reduced Steenrod $\bL$-homology. 
\end{prop}
\begin{proof}
\mk\ni
This sequence can be obtained by adjusting the bounded surgery
theory of \cite{FP} to the continuously controlled case. It is
presented in section 1 of \cite{Pe} in a form where the homology terms are
Ranicki-Wall $L$-groups of the continuously controlled additive
category $\CB(\bar X,Y;\bZ)$. Theorem 2.4 of \cite{Pe} states that
these terms are in fact the reduced Steenrod $\bL$-homology groups of the
corona.
\mk\ni
The naturality follows from the definition of the continuously
controlled category .
\end{proof}

\begin{defn} A subset $X$ of a manifold $M$ has {\it property} $UV^{1}$ if for every neighborhood $U$ of $X$ there is a neighborhood $V$ of $X$ contained in $U$ so that $\pi_{1}(V) \to \pi_{1}(U)$ is trivial. A map $f:M \to Z$ is said to be $UV^{1}$ if each point-inverse $f^{-1}(z)$ is nonempty and $UV^{1}$ in $M$. See section 2 of \cite {La3} for details. 
\end{defn}
\mk\ni
Let $M$ be a closed simply connected $n$-manifold and let $q:M\to Y$ be a $UV^1$-map. Then the
mapping cone $C_q$ is a compactification of $M\times\bR$
by $Y_+=Y\sqcup pt$ which is locally 1-connected at $Y_+$.
Since $(C_q-Y_+)$ is homotopy equivalent to $M$ and $\bar H_*(Y_+;\bL)=H_*(Y;\bL)$,
the controlled surgery exact sequence becomes the following
$$
\dots\to H_{n+1}(Y;\bL)\to\CS^{cc}(C_q,Y_+)\to H_n(M;\bL_{\bullet})\to
H_{n}(Y;\bL)
$$
where $\bL_{\bullet}$ is the connected cover of the spectrum $\bL$. Note that
$(\bL_{\bullet})_0=\GTop$ and by Poincar\'e duality over $\bL_{\bullet}$  \cite{Ra},
$[M,\GTop]=H^0(M,\bL_{\bullet})=H_n(M,\bL_{\bullet}).$
Thus the $n^{th}$ homotopy group
$\CS_n^{cc}(C_q,Y_+)$ of the fiber of the
controlled assembly map of spectra $\bH_*(M;\bL)\to\bH_*(Y;\bL)$ differs from $\CS^{cc}(C_q,Y_+)$ by at most a copy of $\bZ$.
\vskip .2in
\begin{prop}\label{ccforget}
Let $M$ be a simply connected $n$-manifold, $n \ge 5$, and let $q:M\to X$ be a $UV^1$ map and let $M_q$ be its mapping cylinder. Then there is a commutative diagram:
$$
\xymatrix{ \CS^{cc}(C_q,X_+)\ar[rr]^{split- mono}\ar[d]^{forget} &&
\CS_{n+1}^{cc}(C_q,X_+)\ar[r]^{\eta |} &
H_{n+1}(M_q,M;\bL)\ar[d]^{\bar\p}\\
\CS(M)\ar[rr]^{split- mono} &&\CS_n(M) \ar[r]^{\bar\eta} & \bar H_n(M;\bL)\\
}
$$
where $\eta |$ and $\bar\eta$ are isomorphisms.
\end{prop}
\begin{proof}
We have two vertical fibration sequences of spectra on the right, leading to the diagram below at level $n$:
$$
\xymatrix{
&\ar[d]&\ar[d]&\ar[d]\\
&\CS^{cc}(C_{q},X_{+})\ar[d] \ar[r]& \CS^{cc}_{n+1}(C_{q},X_{+}) \ar[d]^{\eta}\ar[r]^{\eta|}& H_{n+1}(M_{q}, M; \bL)\ar[d]^{\bar\p}\\
**[r] H_{n}(M;\bL_{\bullet})\cong&[M\x(0,\,1),\,\GTop]\ar[d]^{q_{*}}\ar[r]&H_{n}(M;\bL)\ar[d]^{q_{*}}\ar[r]^(.5){\cong}&H_{n}(M;\bL)\ar[d]^{q_{*}}\\
&H_{n}(X;\bL)\ar[r]^(.5){\cong}&H_{n}(X;\bL)\ar[r]^(.5){\cong}\ar[d]&H_{n}(X;\bL)\ar[d]\\
&&&
}
$$
from which we see that $\eta|$ is an equivalence. Using the exact sequence 
$$
\xymatrix{
0\ar[r] & H_{n}(M;\bL_{\bullet}) \ar[r]& H_{n}(M;\bL)\ar[r]
&H_{n}(M;\bL/\bL_{\bullet}) \cong\bZ
}
$$
\mk\ni
one chases the diagram above to show first that the composition $\eta':\CS^{cc}(C_{q},X_{+}) \to H_{n+1}(M_{q},M;\bL)$ is a monomorphism and then that its cokernel is a subgroup of $\bZ$. That $H_{n}(M;\bL/\bL_{\bullet}) \cong\bZ$ follows from the observation that $H_{n}(M;\bL/\bL_{\bullet})$ is rationally isomorphic to $\bQ$ and that the only nonzero term in degree $n$ on the $E_{2}$ page of the Atiyah-Hirzebruch spectral sequence computing $H_{n}(M;\bL/\bL_{\bullet})$ is isomorphic to $\bZ$.
\mk\ni
If $[g]\in \CS^{cc}(C_q,X_+)$ is a structure, then there is a manifold $N$ compactified by $X_{+}$ so that $g: N \to M \x (0,\,1)$ extends over $X_{+}$ by the identity and such that this extended map is a strict homotopy equivalence rel $X_{+}$. 
\mk\ni
If $q:M \to X$ is cell-like, then $X$ is locally $k$-connected for all $k$, so there is a retraction from a neighborhood of $X_{+}$ in the compactification of $N$ to $X_{+}$. The proof of the existence of mapping cylinder neighborhoods in \cite{Q} now shows that $X_{+}$ has a mapping cylinder neighborhood in the compactification of N. This splits $N$ as $N' \x (0,\,1)$ and gives a  homotopy equivalence $N' \to M$. The thin $h$-cobordism theorem guarantees that this gives a well-defined forgetful map from $\CS^{cc}(C_{q},X_{+}) \to \CS(M)$. The mapping cylinder projection provides a cell-like map $q_{N}:N \to X$.
\mk\ni
We put $\bar H_n(M;\bL)$, rather than $H_n(M;\bL)$ in the lower right corner of this diagram because both the vertical and horizontal maps have images in $\bar H_n(M;\bL)$ and, as observed in  Proposition \ref{iso}, $\bar \eta$ is an isomorphism.
\end{proof}
\mk\ni
The argument above gives us an important corollary.

\begin{cor}\label{cestructure}
Let $q:M\to X$ be a cell-like map of a simply connected closed manifold $M$.
Then
\begin{enumerate}
\item[(1)]  $\CS^{cc}(C_q,X_+)$ is generated by
strict maps $f:(C_p,X)\to (C_q,X)$ where $p:N\to X$ is a cell-like map.
\item[(2)] The forget control map takes $f$
to a homotopy equivalence $h:N\to M$ which factors through the cell-like maps
$q$ and $p$.
\end{enumerate}
\end{cor}\qed

\ni{\it Proof of Theorem~\ref{main}.}
    (\( \T^{odd}(\CS(M))\subset\CS^{CE}(M).\))\\
\ni    Let $\alpha$ be an odd torsion element of \( \CS(M)\).
 Let $\gamma=\eta(\alpha)\in \bar H_n(M;\bL)$.
By Theorem \ref{1-connected} there is a
cell-like map $q:M\to X$ such that $q_*(\gamma)=0$. Consider the diagram of Proposition \ref{ccforget}.
 By Lemma~\ref{existceappl} there is
a torsion element $\widehat\gamma\in H_{n+1}(M_q,M;\bL)$ such that
$\bar\p(\widehat\gamma)=\gamma$. Let $\alpha'=\eta'^{-1}(\widehat\gamma)$. 
Since $\alpha$ is the image of $\alpha'$ under the forgetful map,
by Proposition  \ref{cestructure} we have $\alpha\in\CS^{CE}(M)$.
\mk
\ni (\( \T^{odd}(\CS(M))\supset\CS^{CE}(M).\))\\
Suppose that \( c:N \to X \) and \( q:M
    \to X \) are cell-like maps and that \( f:N \to M \) is a homotopy
    equivalence such that \( q \circ f \simeq c \).
    \[ \xymatrix{ N \ar[rr]^{f}\ar[dr]^{CE}_{c}&&M\ar[dl]_{CE}^{q}
    \\ & X & }\]
We consider the diagram of Proposition \ref{ccforget}:

$$
\xymatrix{ \CS^{cc}(C_q,X_+)\ar[r]^{\eta'}\ar[d]^{forget} 
 &
H_{n+1}(M_q,M;\bL)\ar[d]^{\bar\p}\\
\CS(M)\ar[r]^{\eta} & \bar H_n(M;\bL)\ar[r] & \bZ\\
}
$$
By the Vietoris-Begle theorem a cell-like map induces an isomorphism
    of ordinary cohomology or Steenrod homology with any coefficients
    (see Proposition \ref{Wa}).
    Therefore \( H_{n}(M;\bL) \to H_{n}(X;\bL) \) is an
    isomorphism rationally and the image of \( H_{n+1}(M_q,M;\bL) \)
    in $ H_{n}(M;\bL)$ is therefore a torsion group.  Since $\bL$ is
    an Eilenberg-MacLane spectrum at 2, \( H_{n}(M;\bL) \to H_{n}(X;\bL) \)
    is an isomorphism at 2 and hence the image of \( H_{n+1}(M_q,M;\bL) \)
    in $ H_{n}(M;\bL)$ is odd torsion.
By Proposition \ref{cestructure} $[f]$ is the image of $[c]\in \CS^{cc}(C_q,
X_+)$ under the forgetful map. Then $[f]=\eta^{-1}\bar\p(\gamma)$ is an odd
torsion element where $\gamma=\bar s([c])\in  H_{n+1}(M_q,M;\bL)$.
\qed

\

\section{Continuous control near the corona}

\mk \ni We move on to the nonsimply connected case. We will use germs of continuously controlled structures near infinity to recover the main diagram in the proof of Proposition \ref{ccforget}. The computation of $\CS^{CE}(M)$ in the nonsimply connected case is made more complicated because we no longer have the isomorphism $\CS(M) \cong \bar H_{n}(M;\,\bL)$. We note, that by Proposition \ref{Wa}, if $q:M \to X$ is cell-like, then there are maps $\xymatrix{ M \ar[r]^{q}& X \ar[r] & E_{2}(M)}$ such that the composition is the inclusion, where $E_{2}(M)$ is the second stage of the Postnikov system of $M$. Elements of $H_{n}(M;\,\bL)$ which survive to $H_{n}(E_{2}(M);\,\bL)$ therefore cannot be in the kernel of $H_{n}(M;\,\bL) \to H_{n}(X;\,\bL)$, so we are led to examine the boundary map $H_{n+1}(E_{2}(M),M;\,\bL) \to H_{n}(M;\,\bL)$, leading to a proof of our main result.

\begin{prop}\label{delta}
Let $(P,Q)$ be a CW pair with an inclusion isomorphism
$\pi_1(Q)=\pi_1(P)=\pi$. Then the homomorphism
$\p':H_{n+1}(P,Q;\bL)\to\CS_n(Q)$ defined in \S2 coincides with
the induced homomorphism on homotopy groups of in the diagram below:
$$
\xymatrix{ \bH_{*+1}(P,Q;\,\bL)\ar[r]\ar[d]^{\p'}&\bH_*(Q;\bL) \ar[r]\ar[d] &\bH_*(P;\bL)\ar[d]^{A_P}\\
           \CS_{*}(Q)\ar[r]&\bH_*(Q;\bL) \ar[r]^{A_{Q}} &\bL_*(\bZ\pi)\\
       }
$$
where $A_P$ and $A_{Q}$  are the assembly maps for $P$ and $Q$.
\end{prop}
\begin{proof} The proof is a diagram chase. \end{proof}
\ni We recall the notation $\delta=p\circ\p'$ where $p:\CS_n(M)\to\CS(M)$
is the projection.
\mk\ni
To prove Theorem \ref{main2} we need the germ version of continuously controlled
surgery theory constructed in \cite{Pe}.

\begin{defn} \
    \begin{enumerate}
            \item [(i)] Let $N$ be an open manifold and let
        $X$ be a compact corona of a compactification
        $\bar N$ of an end of $N$. A {\it strict homotopy equivalence near $X$} is
        a strict map $\bar f:\bar W\to \bar N$, where $\bar W$ is a compactification of an end of $W$ by $X$ and $\bar f|_{ X} = id_{X}$, such that there are neighborhoods $\bar U \supset \bar V$ of $X$ in $\bar N$ and $\bar U' \supset \bar V'$ of $X$ in $\bar W$ such that $f(\bar U) \subset \bar U'$ and there is a strict map $\bar g:\bar U' \to \bar U$ with $\bar g|_{ X} = id_{X}$ such that
        \begin{enumerate}
        \item [(a)] $\bar g \circ \bar f|_{\bar V}$ is strict homotopic in $\bar U$ to $id_{\bar V}.$
        \item [(b)] $\bar f \circ \bar g |_{V'}$ is strict homotopic in $\bar U'$ to $id_{\bar V'}$.
        \end{enumerate}
         \item [(ii)]  Two strict homotopy
        equivalences near $X$, $\bar f:\bar W\to \bar N$
        and $\bar f':\bar W'\to \bar N$ are
        {\it equivalent} if there exist a neighborhood $\bar V$ of $X$ in
        $\bar W$
        and a strict map $h:\bar V \to \bar W'$, $h|_{X}= id_{X}$ which is
        an open  imbedding
        and  $\bar f'\circ h:\bar V\to \bar N$ is  strict homotopic to
        $\bar f|_{\bar V}$.
        \item [(iii)]  The set of the equivalence classes
        of strict homotopy equivalences of manifolds near $X$ is called
        the set of {\it germs of continuously controlled structures}
        on $N$ at $X$ and it is denoted as
        $\CS^{cc}(\bar N,X)_{\infty}$.
    \end{enumerate}
\end{defn}
\ni One can define germs of homotopy classes $[N,\GTop]_{\infty}$ of maps
at $X$ and the corresponding $L$-groups and form a surgery exact sequence.
This was done in \S 15 of \cite{FP} in the case of bounded control and in  \cite{Pe} for continuous control.
We state the result here in the case when $N=M\times(0,1)$ and $\bar N$
is an open mapping cylinder\footnote{This is the usual mapping cylinder of $q$ with the domain copy of $M$ stripped off.}
$\stackrel\circ M_q$ of a cell-like map $q:M\to X$ of a closed orientable
manifold.
\begin{prop}\label{germexact}
Let $q:M\to X$ be a cell-like map of a closed orientable $n$-manifold,
then there is an exact sequence
$$
\dots\to\bar H_{n+1}(X;\bL)\to\CS^{cc}(\stackrel\circ M_q,X)_{\infty}\to
[M,\GTop]\to H_n(X;\bL).
$$
\end{prop}

\ni By Proposition \ref{cestructure}, forget control defines
a map $\phi:\CS^{cc}(\stackrel\circ M_q,X)_{\infty}\to \CS(M)$.
Moreover, there is a commutative diagram:
$$
\xymatrix{
\CS^{cc}(\stackrel\circ M_q,X)_{\infty}\ar[r]\ar[d]^{\phi} & [M,\GTop]\ar[r]\ar[d] &
H_n(X;\bL)\ar[d]^A\\
\CS(M)\ar[r] &[M,\GTop]\ar[r] &L_n(\bZ\pi_1(M)).\\
}
$$
Here $A$ is the assembly map for $X$.

\begin{prop}\label{germforget}
If $q:M\to X$ is a cell-like map of a closed $n$-manifold,
then the forget control map
$\phi:\CS^{cc}(\stackrel\circ M_q,X)_{\infty}\to \CS(M)$
factors as
$$
\CS^{cc}(\stackrel\circ M_q,X)_{\infty}\stackrel{j}\to
H_{n+1}(M_q,M;\bL)\stackrel{i_*}\to
H_{n+1}(E_2(M),M;\bL)\stackrel{\delta}\to\CS(M)
$$
where $j$ is a monomorphism with cokernel $\bZ$ or $0$.
\end{prop}
\begin{proof}
Proposition \ref{Wa} defines a map $g:X\to E_2(M)$ such that
$g\circ q$ is homotopic to the inclusion $M\to E_2(M)$.
We consider the diagram of (horizontal) fibrations of spectra
$$
\xymatrix{ \bH_{*+1}(X,M;\bL)\ar[d]\ar[r] & \bH_*(M;\bL) \ar[rr]\ar[d]^{=} &&\bH_*(X;\bL)\ar[d]^{g_*}\\
\bH_{*+1}(E_{2}(M),M;\bL)\ar[d]\ar[r]&\bH_*(M;\bL) \ar[rr]\ar[d]^{=} &&\bH_*(E_2(M);\bL)\ar[d]^{A_E}\\
 \CS_{*}(M)\ar[r]         &\bH_*(M;\bL) \ar[rr] &&\bL_*(\bZ\pi).\\
       }
$$
In dimension $n$ the homomorphism between homotopy groups of
the fibers gives
$$
H_{n+1}(M_q,M;\bL)\stackrel{i_*}\to H_{n+1}(E_2(M),M;\bL)
\stackrel{\p'}\to \CS_n(M)
$$
where $H_{n+1}(M_q,M;\bL)$ differs from $\CS^{cc}(\stackrel\circ
M_q,X)_{\infty}$ by a potential summand $\bZ$. The proof of this is similar to the proof of Proposition \ref{ccforget}, using the fibration sequence 
$$\CS^{cc}_{*+1}(\stackrel\circ
M_q,X)_{\infty}\to \bH_{n}(M;\,\bL)\to \bH_{n}(X;\,\bL)$$
in place of 
$$\CS^{cc}_{*+1}(C_{q},X_{+})\to \bH_{n}(M;\,\bL)\to \bH_{n}(X;\,\bL)$$
The result then follows from Proposition \ref{delta}.
\end{proof}

\ni{\it Proof of Theorem \ref{main2}.} (\( \CS^{CE}(M)\supset im(\delta^T_{[2]}).) \) We are given an odd torsion element \( \alpha\in H_{n+1}(E_2(M),M; \bL) \) with \(\delta(\alpha)= [f] \in \CS(M) \) where
     \(\delta\) is the composition
    \[ H_{n+1}(E_2(M),M; \bL) \cong \mathcal{S}_{n+1}(E_2(M),M)\to
    \mathcal{S}_n(M)\to\CS(M). \]
    \ni By
    Proposition \ref{existceappl}, there exist a cell-like map \( q:M \to X \)
    and an odd torsion element \( \widehat \alpha \in H_{n+1}(M_q,M;\bL)\cong \CS_{n+1}^{cc}(\stackrel\circ M_q,X)_{\infty} \) so that \( \alpha
    \) is the image of \( \widehat \alpha \) under the inclusion-induced
    map \(i_*: H_{n+1}(M_q, M; \bL ) \to H_{n+1}(E_2(M), M; \bL)\).
    Since $\widehat \alpha$ has finite order and $j:\CS^{cc}(\stackrel\circ M_q,X)_{\infty} \to \CS_{n+1}^{cc}(\stackrel\circ M_q,X)_{\infty}$ is an isomorphism on torsion subgroups, $ \widehat\alpha=j(\alpha')$, where $\alpha' \in \CS^{cc}(\stackrel\circ M_q,X)_{\infty}$. By
    Propostion \ref{germforget} $\phi(\alpha')=[f]$.
    Let $g:W\to M\times(0,1)$ be a representative for $\alpha'$.
    As in the proof of Corollary \ref{cestructure} we may assume that
    $W=N\times(0,1)$ and $\bar W=M_p$, where $p:N\to X$ is cell-like.
    Thus, $[f]=\phi(\alpha')$  is realized by cell-like maps $p$ and $q$.
\vskip .2in

\ni(\( \CS^{CE}(M)\subset im(\delta^T_{[2]})). \) Suppose that \( c:N \to X \) and \( q:M
    \to X \) are cell-like maps and that \( f:N \to M \) is a homotopy
    equivalence such that \( q \circ f \simeq c \).
    \[ \xymatrix{ N \ar[rr]^{f}\ar[dr]^{CE}_{c}&&M\ar[dl]_{CE}^{q}
    \\ & X & }\]
    \ni As above,
    there is an inclusion-induced map \( p:X \to E_2(M) \) and the forgetful map \(
    H_{n+1}(M_q,M;\bL) \cong\CS_{\infty}^{cc}(M_q,X) \to \CS(M) \) factors through \(
    H_{n+1}(E_2(M),M;\bL) \).  It therefore suffices to show that the image
    of \( H_{n+1}(M_q,M;\bL)\)  in \( H_{n+1}(E_2(M),M;\bL) \) is an odd torsion
    group.  By the Vietoris-Begle theorem a cell-like map induces an isomorphism
    of ordinary cohomology or Steenrod homology with any coefficients
    (see Proposition \ref{Wa}).
    Therefore \( H_{*}(M;\bL) \to H_{*}(X;\bL) \) is an
    isomorphism rationally, and hence, the image of \( H_{*}(M_q,M;\bL) \)
    in  \( H_{*}(E_2(M),M;\bL) \) is torsion.  Since $\bL$ is
    an Eilenberg-MacLane spectrum at 2, \( H_{*}(M;\bL) \to H_{*}(X;\bL) \)
    is an isomorphism at 2 and hence \( H_{*}(E_2(M),M;\bL) \) is odd torsion.
\qed

\section{Deforming Riemannian manifolds in Gromov-Hausdorff space}

\mk\ni In this section, we apply the theory developed above to study Riemannian manifolds in Gromov-Hausdorff space.

\ni\begin{defn}\label{contr}\

\begin{enumerate}

        \item [(i)] A continuous function $\rho:\bR_+\to\bR_+$
with $\rho(0)=0$, continuous at $0$, with $\rho(t)\ge t$ for all $t$ is a {\it contractibility
function} for a metric space $X$ if
there is $R>0$ such that for each $x\in X$ and $t\le R$, the $t$-ball $B_t(x)$
centered at $x$ can be contracted to a point in the $\rho(t)$-ball
$B_{\rho(t)}(x)$. 
        \item [(ii)] Similarly, if $X$ is an n-manifold, $\rho$ is a {\it topological injectivity function} for $X$ if for each $x \in X$ and $t \le R$ there is an open subset $U \subset X$ so that $U$ is homeomorphic to $\bR^{n}$ and $B_{t}(x) \subset U \subset B_{\rho(t)}(x)$.
    \end{enumerate}
\end{defn}

\mk\ni Let $\rho=\rho_{1}:[0,\,R) \to [0,\,\infty)$ be a contractibility function. The theorem of Petersen on page 392 of \cite{Pet} shows that for every $\epsilon>0$ there is a $\delta>0$ so that if $X$ and $Y$ are compact $n$-dimensional metric spaces with contractibility function $\rho$ such that $d_{GH}(X,\,Y)<\delta$, then $X$ and $Y$ are homotopy equivalent by maps and homotopies that move points by less than $\epsilon$. Moreover, given $\epsilon$, there is an explicit computation of the necessary $\delta$.
\mk\ni
Combining this with the results of Chapman-Ferry, Freedman-Quinn, and Perelman cited in the introduction, we see that if $M$ is a closed topological $n$-manifold with a given metric $d_{M}$, and a contractibility function $\rho$,  then there is a $\delta>0$ such that any other topological $n$-manifold with contractibility function $\rho$ and $d_{GH}(M,N)<\delta$ must be homeomorphic to $M$.
\mk\ni
In this section, we show that the condition that $M$ be stationary with a fixed metric is necessary: that there are families of nonhomeomorphic Riemannian manifolds with a common contractibility function that can be deformed arbitrarily close to each other in a precompact region of Gromov-Hausdorff space. We get a complete algebraic description of this behavior and produce many examples of nonhomeomorphic families of manifolds with common topological injectivity functions that can be similarly pushed together.

\newpage
\begin{defn}\label{GH} \

    \begin{enumerate}
    	 \item [(i)] If $Z$ is a metric space, $X \subset Z$, and $\epsilon > 0$, 
	 $N_{\epsilon}(X)=\{ z \in Z\,|\,d(z,X)<\epsilon\}$.
        \item [(ii)] If $X$ and $Y$ are compact subsets of a metric space
        $Z$, the {\it Hausdorff distance} between $X$ and $Y$ is
        $$ d_H(X,Y)=\inf\{\epsilon>0\mid X\subset N_{\epsilon}(Y),
        Y\subset N_{\epsilon}(X)\}.$$
        Here, $X$ and $Y$ are isometrically embedded in $Z$.
        \item [(iii)] If $X$ and $Y$ are compact metric spaces, the
        {\it Gromov-Hausdorff distance} from $X$ to $Y$ is
        $$ d_{GH}(X,Y)=\inf_Z\{d_H(X,Y)\mid X,Y\subset Z\}.$$
        \item[(iv)] Let $\cm$ be the set of isometry classes of compact
        metric spaces with the Gromov-Hausdorff metric.
        \item[(v)] Let $\mn$ be the set of all $(X,d)\in\cm$ such that
        $X$ is a topological $n$-manifold with (topological) metric $d$
        with contractibility function $\rho$.
        \end{enumerate}
\end{defn}
\ni It is well-known that $\cm$ is a complete metric space (see \cite{Gro} or
\cite{Pet} for an exposition).
\begin{thm}\label{push}\
\begin{enumerate}
\item[(i)] If $n\ne 3$ and $X\in\cm$ is in the closure of $\mn$, then there is
an $\epsilon>0$ so that there are only finitely many homeomorphism types of
manifolds $M\in\mn$ with $d_GH(M,X)<\epsilon$. If $d_{GH}(M,X)$,
$d_{GH}(M',X)<\epsilon$, then there exists a simple homotopy equivalence
$h:M'\to M$ which preserves rational Pontryagin classes.
\item[(ii)] If $[f] \in \CS^{CE}(M)$ with $f:N \to M$, $M$ and $N$ smooth, then 
there exist a contractibility function $\rho$ and a compact metric space $X$ such that  every
neighborhood of $X$ in $\cm$ contains manifolds lying in $\mn$
and homeomorphic to both $M$ and $N$.
\item[(iii)] There exist examples as in (ii) such that $M$ and $N$ are not homeomorphic.

\end{enumerate}
\end{thm}
\begin{proof}
Part (i) is Theorem 2.10 of \cite{F1}.
\mk\ni
For part (ii),  and let $q:M\to X$ and
$p:N\to X$ be cell-like maps. By the main results of \cite{FO} and \cite{M},
there exist a contractibility function $\rho$, and sequences of
Riemannian metrics $\{d_i^M\}$ and $\{d_i^N\}$ on $M$ and $N$
respectively lying in $\mn$ and converging in $\cm$ to $(X,d)$ for
some metric $d$.
\mk\ni
For part (iii), let $M$ and $N$ be the manifolds from Corollary \ref{two} or Proposition \ref{Diarmuid}.
\end{proof}
\mk\ni
Let $\overline{\mn}$ be the closure of $\mn$ in Gromov-Hausdorff 
space and let $\p\mn$ be the boundary. 
For a compact metric space $X$, we will denote its
isometry class by the same letter $X$.

\begin{thm}\label{converse}
Suppose that the isometry type of a metric space $X$ belongs to
$\p\mn$. Then there is an $\epsilon>0$ such every two manifolds
$M,N\in B_{\epsilon}(X)\cap \mn$ are CE-related.
\end{thm}
\begin{proof} The proof will follow Proposition \ref{tensor}.\end{proof}
\begin{defn}
A map $f:M\to X$ has the {\it $\delta$-lifting property in dimensions
$\le k$} if for every PL pair $(P,Q)$, $\dim P\le k$ for every
commutative diagram
\[
    \xymatrix{ Q \ar[rr]^{g'}\ar[d]& &M\ar[d]^{f}\\
    P \ar[rr]^{g}\ar[urr]^{\bar g} & & X. }
    \]
there is a map $\bar g:P\to M$ extending $g'$ such that
$\dist(f\circ\bar g,g)<\delta$.
\end{defn}
\begin{prop}\label{hisom}
Let $X$ be a locally $k$-connected space for $k>n$, then there
exists $\delta>0$ such that every map $f:Z\to X$ from a compact
$n$-dimensional ANR with the $\delta$-lifting property in
dimensions $\le n+1$ is a weak homotopy equivalence through dimension
$n$ (i.e., such that $f$ is $n+1$-connected). Furthermore, $f$ induces
isomorphisms of Steenrod homology groups $f_*:H_i(M)\to
H_i(X)$ for $i\le n$.
\end{prop}
\begin{proof}
The weak homotopy equivalence in dimension $n$ easy follows from
the lifting property. This implies the result for singular
homology. We note that the Steenrod homologies coincide with the
singular homologies in the locally $k$-connected case.
\end{proof}
\begin{prop}\label{map}
If $X\in\p\mn$, then for every $\delta>0$ there exists
$\epsilon>0$ such that for every $M\in\mn$ with $d_{GH}(M,X)<\epsilon$
there is a map $f:M\to X$ with the $\delta$-lifting property in
dimensions $\le n+1$.
\end{prop}
\begin{proof}
The space $X$ is locally $k$-connected for all finite $k$ (see
\cite{F3}). Then for small $\epsilon$ a map $f:M\to X$ can be
constructed by induction by means of a small triangulation on $M$
(if $M$ does not admit a triangulation, one can use a CW complex
structure). Given $\delta_0>0$, we may assume that
$d(x,f(x))<\delta_0$. Clearly, for a proper choice of $\delta_0$
the map $f$ will have the $\delta$-lifting property.
\end{proof}
\mk\ni
We refer to Madsen and Milgram~\cite{MM} for the following:
\begin{prop}\label{MM}
For any finite complex $K$ there is an isomorphism
$$
H_n(K;\bL_{(2)})\cong\bigoplus_iH_{n+4i}(K;\bZ_{(2)})\oplus H_{n+4i-2}(K;\bZ_2)
$$
which is natural with respect to maps $K\to L$.
\end{prop}
\begin{prop}\label{LQ-isom}
If $X\in\p\mn$, then there exists $\epsilon>0$ such that for every
$M\in\mn$ with $d_{GH}(M,X)<\epsilon$ there is a map $f:M\to X$
such that $f_*:H_*(M;\bL_{(2)}) \to
H_*(X;\bL_{(2)})$ is an isomorphism.
\end{prop}
\begin{proof}
Since $\bL_{(2)}$ is an Eilenberg-MacLane spectrum, 
we can take $\epsilon$ from Proposition \ref{map}. Then  Proposition
\ref{hisom} and the fact that $H_i(M)=H_i(X)=0$ for $i>n$ imply
the required result. This last fact follows from the arguments above. Given 
$k$, one shows that there is a $\delta > 0$ so that if $d_{GH}(M^{n},X)<\delta$, there 
is a $k$-connected map $M \to X$. Repeating this for a sequence of $M$'s shows that the 
 homology groups of $X$ are trivial in dimensions $>n$. \end{proof}
\mk\ni
Petersen \cite{Pet} correctly concludes from similar arguments that $X$ can have no finite-dimensional subsets of dimension $>n$ and incorrectly concludes from 
 this that $X$ must have covering dimension $\le n$. The limit spaces $X$ constructed in this paper are infinite-dimensional spaces containing no finite-dimensional subspaces of dimension $>n$. See \cite{Wa2}, \cite{Dr} for further explanation.
\mk\ni
The following is well-known:
\begin{prop}\label{tensor}
Suppose that a homomorphism of abelian groups $\phi:A\to B$ becomes trivial after tensoring with $\bZ_{(2)}$. Then
the image of $\phi$ lies in the odd torsion subgroup of $B$, $\phi(A)\subset T^{odd}(B)$.
\end{prop}

\ni{\it Proof of Theorem \ref{converse}.} We take $\epsilon$ from Proposition \ref{LQ-isom}. Let $c:N\to X$
and $q:M\to X$ be corresponding maps. We may assume that there is
a homotopy lift $f:N\to M$ of $c$ which is a homotopy equivalence. Then $f$ induces isomorphisms $f_*:H_*(N;\bL_{(2)}) \to
H_*(M;\bL_{(2)})$.
\mk\ni
As was shown in \cite{F3} (P4, page 98), there are finite polyhedra $P_1$, $P_2$ and maps $p_1:X\to P_1$,
$p_2:X\to P_2$ and $g:P_2\to P_1$ such that $p_1=g\circ p_2$, $p_2$ is 
$n+3$-connected and $g$ is $(\dim P_1+3)$-connected. Let $q_i=p_i\circ q$, $i=1,2$.
We note that these conditions imply that $q_2$   induces isomorphisms of homology in dimension $\le n+3$ and
$g$ induces isomorphisms of homology in dimension $\le \dim P_1+3$. The latter implies that $\im g_*=\im (q_1)_*$
for homology of dimension $\le \dim P_1+3$. In view of Proposition~\ref{MM} we obtain the following:
\begin{enumerate}
\item \ \  $(q_2)_*:H_n(M;\bL_{(2)})\to H_n(P_2;\bL_{(2)})$ is a monomorphism;
\item \ \  $\im g_*=\im (q_1)_*$ for the $(n+1)$-dimensional $\bL_{(2)}$-homology.
\end{enumerate}
\mk\ni
We claim that $g_*:H_{n+1}(P_2,M;\bL_{(2)})\to H_{n+1}(P_1,M;\bL_{(2)})$ is the zero homomorphism.
\mk\ni
Consider the commutative diagram generated by exact sequences of pairs
$$
\begin{CD}
H_{n+1}(M;\bL_{(2)}) @>(q_2)_*>> H_{n+1}(P_2;\bL_{(2)}) @>j^2_*>>  H_{n+1}(P_2,M;\bL_{(2)}) @>\partial_2>>\\
@V=VV @Vg_*VV  @Vg_*VV\\
H_{n+1}(M;\bL_{(2)}) @>(q_1)_*>> H_{n+1}(P_1;\bL_{(2)}) @>j^1_*>>  H_{n+1}(P_1,M;\bL_{(2)}) @>\partial_1>>\ .\\
\end{CD}
$$
Let $\alpha\in H_{n+1}(P_2,M;\bL_{(2)})$. By the property (1) $\partial_2(\alpha)=0$. Hence $\alpha=j^2_*(\beta)$
for some $\beta$. By the property (2) there is $\gamma\in H_{n+1}(M;\bL_{(2)})$ such that $(q_1)_*(\gamma)=g_*(\beta)$.
Hence $0=j_*^1\circ g_*(\beta)=g_*(\alpha)$ and the claim is proven.
\mk\ni
Since $H_*(Y;\bL_{(2)})=H_*(Y;\bL)\otimes\bZ_{(2)}$, Proposition~\ref{tensor} implies that $g_*$ takes 
$H_{n+1}(P_2,M;\bL)$ to odd torsion.
\mk\ni
Since $p_i\circ
q$ is 2-connected, the space $E_2(M)$ can be constructed out of
$P_i$ by killing higher dimensional homotopy groups. Thus the
inclusion $M\subset E_2(M)$ can be factored through $X$ and $P_i$, $i=1,2$. Hence there is a commutative diagram
\begin{equation}
\xymatrix{
    H_{n+1}(P_{2},M;\bL) \ar[rr]^{\cong}
    \ar[d] & &\CS_{n+1}(P_2,M)\ar[rr]^{\phi_2}\ar[d] & & \CS_n(M)\ar[d]^{=}\\
H_{n+1}(P_{1},M;\bL) \ar[rr]^{\cong}
    \ar[d] & &\CS_{n+1}(P_1,M)\ar[rr]^{\phi_1}\ar[d] & & \CS_n(M)\ar[d]^{=}\\
H_{n+1}(E_2(M),M;\bL) \ar[rr]^{\cong} & & \CS_{n+1}(E_2(M),M)\ar[rr]^{\partial}
    & & \CS_n(M) . }
\end{equation}
\mk\ni
We show that the element $[f]\in\CS(M)$ comes from an odd torsion element of $H_{n+1}(E_2(M),M;\L)$. 
By Theorem 2.6 of \cite{F3} the structure $[f]$ defined by
$f:N\to M$ belongs to
the kernel of the induced map $(q_2)_*:\CS_n(M)\to
\CS_n(P_2)$. Thus $[f]\in \im(\phi_i)$, $i=1,2$. By the above $\phi_2$ factors through odd torsion.
Therefore $[f]$ is the image under $\partial$ of an odd torsion element. Hence, $[f]\in \im(\delta^T_{[2]})$.
\mk\ni
Applying Theorem~\ref{main2} shows that $N$ and $M$ are CE-related.
\qed
\mk\ni
Next, we demonstrate the existence of topological injectivity functions. Our argument is an easy modification of McMillan's Cellularity Criterion \cite{McM}, which says that a compact subset $X$ of a closed manifold $M^{n},\ n\ge 5$ is a nested intersection of open sets homeomorphic to $\bR^{n}$ if and only if $X$ is cell-like and for every open neighborhood $U$ of $X$ there is an open neighborhood $V$ of $X$ contained in $U$ such that the inclusion induced map $\pi_{1}(V-X) \to \pi_{1}(U-X)$ is trivial. \footnote{Using work of Perelman and Freedman, the cellularity criterion is now known to be true in all dimensions.}
\mk\ni
\begin{thm}\label{topinj}
Let $M$ be a closed topological $n$-manifold with a contractibility function $\rho:[0,\,R)\to [0,\,\infty)$. If $Q$ is a closed $k$-manifold, $k \ge 1,\ n+k \ge 5$, then $M \x Q$ has a topological injectivity function.
\end{thm}
\begin{proof}
\ni McMillan shows that if $X$ is cell-like and $M_{i},\,i=0,\ldots,3$ are nested compact PL manifolds, $M_{i+1}\subset M_{i}$, containing $X$ with the inclusion of $M_{i+1}$ into $M_{i}$ nullhomotopic, $i=0,\ldots,2$, and the inclusion induced map $\pi_{1}(\circM_{3}-X) \to \pi_{1}(\circM_{2}-X)$ is zero, then there is an open set $U$ with $M_{3} \subset U \subset M_{0}$ and $U$ homeomorphic to $\bR^{n}$.\footnote{McMillan's argument shows $X \subset U \subset M_{0}$, but adding one more layer gives the stated result.}
\mk\ni
By immersion theory, any open subset of a topological manifold that contracts to a point in that manifold has a PL structure, so for any $x \in M$, the contractibility function allows us to find arbitrarily long nested sequences of compact PL manifold neighborhoods $M_{i}$ of $x$ with $M_{i+1}\to M_{i}$ nullhomotopic for all $i$. A bit of manipulation to fill in holes, see \cite{McM}, allows us to assume that each $M_{i}$ is connected with connected boundary and that $M_{i}-\circM_{i+1}$ is connected for each $i$. If $q \in Q$ is written as a nested intersection $\{q\}=\cap B_{i}$ of balls in $Q$, then the sequence $M_{i}\x B_{i}$ satisfies McMillan's conditions and guarantees the existence of a topological injectivity function. With a bit of care, this function $\tau$ can be written explicitly in terms of $\rho$.
\end{proof}

\begin{rem}
For each $x \in M \x Q$, this allows us to construct a sequence of
homeomorphisms between euclidean neighborhoods of $x$ and euclidean
neighborhoods of nearby points in $M' \x Q$, where $M$ is deformable to
$M'$. Evidently, these homeomorphisms cannot be controlled well enough
to stitch them together to provide an isomorphism of tangent
microbundles, since that would contradict the characteristic class
computations of Corollary \ref{two}. This suggests that there is no
reasonable way of assigning a tangent bundle to the infinite-dimensional
(but finite cohomological dimensional) homology manifold $X$.
Nevertheless, the controlled Mischenko-Ranicki symmetric signatures
$\Delta_{M}\in H_{n}(M;\,\bL)$ and $\Delta_{M'}\in H_{n}(M';\,\bL)$ map
to the same class in $H_{n}(X;\,\bL)$ after inverting 2. This suggests
that $X$ may possess a well-defined characteristic class theory.\end{rem}

\end{document}